\documentclass[10pt,a4paper]{smfart}

\usepackage[latin1]{inputenc}
\usepackage[T1]{fontenc}
\usepackage[francais]{babel}
\usepackage{amsfonts,diagrams}
\usepackage{euscript}
\usepackage{graphicx}
\usepackage[vcentermath]{youngtab}

\newtheorem{conj}{Conjecture}

\newtheorem{theo}{Th\'eor\`eme}
\newtheorem{prop}{Proposition}

\newtheorem{lemm}{Lemme}

\newtheorem{defi}{D\'efinition}

\newtheorem{nota}[defi]{Notation}

\newtheorem{rema}[defi]{Remarque}

\newcommand{\res}[2]{\vspace{.15cm}

\noindent
{\bf #1 :} {\it #2} \vspace{.15cm} 

\noindent \hspace{-.32cm}
}

\begin{document}

\newcommand{\Es}[4]{ {\ }^{#1} \hspace{-2pt} E^{#2,#3}_{#4} }
\newcommand{\Ed}[5]{ {\ }^{#1,#2} \hspace{-2pt} E^{#3,#4}_{#5} }
\newcommand{\scal}[1]{\langle #1 \rangle}
\newcommand{\ssi}{si et seulement si }
\newcommand{\litt}{Littlewood-Richardson}
\newcommand{\as}{${\cal A}_6\ $}
\newcommand{\mt}{${\cal M}_3\ $}
\newcommand{\Pro}{{\bf I \hspace{-3pt} P}}
\newcommand{\p}{{\mathbb P}}
\newcommand{\C}{\mathbb{C}}
\newcommand{\N}{\mathbb{N}}
\newcommand{\NN}{{\cal N}}
\newcommand{\Z}{\mathbb{Z}}
\newcommand{\oo}{{\cal O}}
\newcommand{\ooo}{{\cal O}}
\newcommand{\hlm}{{}^h\lambda^-}
\newcommand{\vlm}{{}^v\lambda^-}
\newcommand{\cent}[1]{\begin{center} \vspace{-.13cm} #1 \end{center}
\vspace{-.13cm}}
\newcommand{\PM}{P_{\mbox{max}}}
\newcommand{\QM}{Q_{\mbox{max}}}
\newcommand{\sgn}{\mbox{sgn}}
\newcommand{\n}{\mathfrak n}
\newcommand{\B}{\mathfrak b}
\newcommand{\nd}{$\frac{n}{2}\,$}
\newcommand{\nq}{$\frac{n}{4}\,$}
\newcommand{\f}{${\cal F}(Q_1)\ $}
\newcommand{\fp}{${\cal F}'(Q_1)\ $}
\newcommand{\fx}{${\cal F}_x(Q_1)\ $}
\newcommand{\fpx}{${\cal F}_x'(Q_1)\ $}
\newcommand{\se}{Sec(X)-X}
\newcommand{\sv}{Severi variety }
\newcommand{\vs}{Severi variety }
\newcommand{\g}{\mathfrak g}
\newcommand{\h}{\mathfrak h}
\newcommand{\cal}[1]{\mathcal{#1}}
\newcommand{\suiteexacte}[3]{
0 \rightarrow #1 \rightarrow #2 \rightarrow #3 \rightarrow 0 }

\newcommand{\lpara}{\vspace{4pt}

\noindent}
\newcommand{\para}{\ \vspace{8pt}

\noindent}
\newcommand{\Para}{\vspace{12pt}

\noindent}
\newcommand{\PPara}{\vspace{25pt}

\noindent}

\newcommand{\produit}[2]{\oplus_{\alpha_i : \sum \alpha_i = #1} \otimes_i
Z^{k_i-\alpha_i-1,k_i} #2}
\newcommand{\dem}{\noindent \underline {\bf D\'emonstration :} }
\newcommand{\rem}{\noindent \underline {\bf Remarque :} }
\newcommand{\fin}{\begin{flushright} \vspace{-16pt}
 $\bullet$ \end{flushright}}

\newcommand{\surmap}{\longrightarrow \hspace{-.5 cm} \longrightarrow}
\newcommand{\pisur}[2]{
\pi_1(#1) \surmap \pi_1(#2)   } 

\newcommand{\fonction}[5]{
\begin{array}{rrcll}
#1 & : & #2 & \rightarrow & #3 \\
   &   & #4 & \mapsto     & #5
\end{array}  }
\newcommand{\fonc}[3]{
#1 : #2 \rightarrow #3 }

\newcommand{\DC}[8]{
\begin{array}{ccc}
#1             & \stackrel{#2}{\longrightarrow} & #3 \\
\downarrow #4  &                                & \downarrow #5 \\
#6             & \stackrel{#7}{\longrightarrow} & #8
\end{array} }

\newcommand{\matdd}[4]{
\left (
\begin{array}{cc}
\mbox{}#1 & \mbox{}#2  \\
\mbox{}#3 & \mbox{}#4
\end{array}
\right )   }

\newcommand{\mattt}[9]{
\left (
\begin{array}{ccc}
\mbox{}#1 & \mbox{}#2 & \mbox{}#3 \\
\mbox{}#4 & \mbox{}#5 & \mbox{}#6 \\
\mbox{}#7 & \mbox{}#8 & \mbox{}#9
\end{array}
\right )   }

\newcommand{\vectq}[4]{
\left (
\begin{array}{r}
#1 \\
#2 \\
#3 \\
#4
\end{array}
\right )   }

\title{Th\'eor\`emes d'annulation et lieux de d\'eg\'en\'erescence en petit corang}
\author{Pierre-Emmanuel Chaput\\
        Pierre-Emmanuel.Chaput@math.univ-nantes.fr}
\maketitle

\begin{center}
\bf{Abstract}
\footnote{
{\it Classification AMS\/}: 14F17.\\
{\it Mots-cl\'es\/}: Th\'eor\`emes d'annulation, grassmannienne.}
\end{center}

After giving an explicit description of all the non vanishing Dolbeault 
cohomology groups of ample
line bundles on grassmannians, I give two series of
vanishing theorems for ample vector
bundles on a smooth projective variety. They imply
a part of a conjecture by Fulton and Lazarsfeld about the connectivity of some
degeneracy loci.

\vspace{1cm}

\begin{center}
\bf{Introduction}
\end{center}

Soit $X$ une vari\'et\'e projective complexe, $E$ un fibr\'e vectoriel sur
$X$ et $L$ un fibr\'e en droites. Supposons que l'on ait une forme quadratique
sur $E$ \`a valeurs dans $L$, soit une section de $S^2E^* \otimes L$. Si $k$
est un entier, on note $D_k(E)$ le sous-sch\'ema de $X$ o\`u 
cette section est au plus de rang $k$.
Dans \cite[Remark 2,p.50]{fulton}, on peut lire la conjecture suivante
($t(x):=\frac{x(x+1)}{2}$):
\begin{conj}
Soit $E$ un fibr\'e vectoriel
de rang $e$, sur une vari\'et\'e $X$ lisse, projective, connexe 
et de dimension $n$. Supposons que $E$ est
muni d'une forme quadratique \`a valeurs dans un fibr\'e en droites $L$. 
Soit $k$ un entier et supposons que
\begin{itemize}
\item
$\dim D_k(E)=\rho:=n-t(e-k) \geq 1$.
\item
$S^2E^* \otimes L$ est ample.
\end{itemize}
Alors, $D_k(E)$ est connexe.
\label{conjecture}
\end{conj}

Dans cet article, je montre cette conjecture sous l'hypoth\`ese suppl\'ementaire
que $e-k \leq 4$ et $\rho \geq 2$. J'obtiens en fait les r\'esultats 
plus pr\'ecis suivants:
\res{Th\'eor\`eme \ref{degenerescence}}
{
Sous les hypoth\`eses de la conjecture pr\'ec\'edante, \`a part que $X$ n'est plus
suppos\'ee connexe, et si de plus $\rho \geq 2$ et $e-k \leq 4$, alors
l'application de restriction 
$H^q(X,{\cal O}_X) \rightarrow H^q(D_k(E),{\cal O}_{D_k(E)})$ est un
isomorphisme pour $0 \leq q<\rho-1$, et est injective pour $q=\rho-1$.
}

La conjecture de Fulton et Larzarsfeld a \'et\'e r\'esolue en utilisant
une technique diff\'erente par Tu \cite{tu1}, dans le cas o\`u $k$ est pair. Par
ailleurs, \cite{tu2} montre
aussi la connexit\'e, mais \`a condition que $\rho \geq e-k$. Si $e-k=3$ ou 4, ma
borne, 2, est donc meilleure. J'esp\`ere aussi que ce travail est un pas de plus
vers une parfaite compr\'ehension des ph\'enom\`enes combinatoires qui permettent
d'\'etablir les th\'eor\`emes d'annulation. Mentionnons enfin que 
\cite[prop 2.3]{laytimi} donne en utilisant la même technique que moi
la connexit\'e d'un lieu de dimension
strictement positive, mais pour $e-k \leq 2$, ce qui, comme le lecteur pourra
le mesurer, simplifie consid\'erablement le probl\`eme.

\lpara

Il est bien connu que la conjecture d\'ecoule de th\'eor\`emes d'annulation ad\'equats
\cite{manivel}: il suffit en effet d'utiliser une r\'esolution 
du faisceau structural de
$D_k(E)$ par des fibr\'es vectoriels qui sont des puissances de Schur de $E$, et
d'appliquer les th\'eor\`emes d'annulation aux termes de cette r\'esolution.
Dans le cas o\`u $e-k \leq 2$, seuls des crochets, \`a savoir des partitions dont
seule la premi\`ere part est \'eventuellement sup\'erieure \`a 2, 
interviennent, et le th\'eor\`eme 2.1 dans \cite{nahm} convient.
Dans cet article, je propose une g\'en\'eralisation de
ce r\'esultat pour des produits tensoriels de crochets. Cette g\'en\'eralisation 
donne bien entendu un th\'eor\`eme d'annulation pour toutes les partitions; 
malheureusement, la borne obtenue est insuffisante pour \'etablir la conjecture.
Je propose donc dans cet article une m\'ethode un peu nouvelle, de ``comparaison
de suites spectrales'', pour \'etablir des th\'eor\`emes d'annulation plus puissants.
Le th\'eor\`eme \ref{degenerescence} est alors cons\'equence des th\'eor\`emes
d'annulation obtenus par cette m\'ethode.

L'efficacit\'e des th\'eor\`emes d'annulation obtenus d\'epend directement de notre
compr\'ehension de la cohomologie des fibr\'es en droites homog\`enes amples sur une
grassmannienne; Snow \cite{snow_classique} a donn\'e pour la calculer 
une m\'ethode 
diagrammatique commode. L'inconv\'enient de cette m\'ethode est qu'il est difficile
d'en d\'eduire, pour $r<e$ et $l$ des entiers fix\'es, l'entier maximal
$p$ tel qu'il existe $q$ avec $H^{p,q}[G(r,e),\ooo(l)] \not = 0$
($G(r,n)$ d\'esigne la grassmannienne des $r$-plans vectoriels dans un espace
vectoriel de dimension $n$ et $\ooo(l)$ est la $l$-i\`eme puissance du 
d\'eterminant du fibr\'e quotient). Pourtant, 
comme nous allons le voir, la d\'etermination de cet entier est cruciale pour
obtenir des th\'eor\`emes d'annulation. Je propose donc une description nouvelle de
la cohomologie de $\ooo(l)$ sur une grassmannienne, bas\'ee sur celle de Snow, et
j'en d\'eduis le th\'eor\`eme \ref{partition_p_max} 
qui permet de d\'eterminer cet entier.

\lpara

Je remercie mon directeur de th\`ese Laurent Manivel pour son aide tout au long 
de l'\'elaboration de cet article.
Cet article a \'et\'e publi\'e \`a Documenta Mathematica \cite{documenta}.

%******************************************************************************

%******************************************************************************

\section{Cohomologie de Dolbeault des fibr\'es $\ooo(l)$ sur une grassmannienne}

%******************************************************************************

\subsection{Description explicite de toutes les partitions admisibles}

Rappelons tout d'abord la description de Snow \cite{snow_classique} des groupes de
cohomologie non nuls sur une grassmannienne. Une partition est suite 
d\'ecroissante finie d'entiers. Je note $\lambda_i$ 
le $i$-i\`eme entier de $\lambda$; c'est par
convention la $i$-i\`eme part de $\lambda$. Le poids de $\lambda$, not\'e
$|\lambda |$, est la somme de ses parts. Sa longueur, $l(\lambda)$, est le
nombre de parts non nulles.
Dans de
nombreuses circonstances, il est pratique de repr\'esenter les partitions
par un diagramme comme le sugg\`ere l'exemple qui suit:
$$
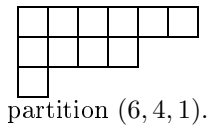
\begin{array}{c}
\yng(6,4,1)\\
\mbox {partition } (6,4,1).
\end{array}
$$
A chaque partition $\lambda$ correspond un foncteur
de la cat\'egorie des espaces vectoriels dans elle-même
que j'appelle ``de Schur'' et que je 
note $S_\lambda$ \cite{fulton_harris}.

Lorsque $\lambda$ est une partition, on peut attribuer \`a chaque case de cette
partition son nombre de crochet qui est le nombre de cases de $\lambda$
situ\'ees en-dessous ou \`a droite de cette case, cette case comprise. Par exemple,
les cases de la partition suivante ont \'et\'e num\'erot\'ees 
par leurs nombres de crochets:
$$
\young(6421,31,1).
$$
Convenons alors, si $l$ est un
entier, que $\lambda$ sera appel\'ee $l$-admissible si toutes les cases
re\c{c}oivent un num\'ero diff\'erent de $l$.
 
Snow \cite{snow_classique} a montr\'e qu'il existe une bijection, 
\`a $r<e$ et $l$ fix\'es, entre
les composantes des groupes de cohomologie $H^{p,q}[G(r,e),\ooo (l)]$
et les partitions de
longueur inf\'erieure ou \'egale \`a $r$ dont toutes les parts sont
inf\'erieures ou \'egales \`a $e-r$ (par la suite je les appellerai simplement
partitions de taille $(e-r,r)$) qui sont $l$-admissibles. 
Si $\lambda$ est une partition $l$-admissible, 
soit $\hlm$ (respectivement $\vlm$)
la suite d'entiers telle que
 $\hlm_i$ (respectivement $\vlm_j$) est le nombre de cases 
situ\'ees 
sur la $i$-i\`eme ligne (respectivement sur la $j$-\`eme colonne) dont le 
nombre de crochet est strictement inf\'erieur \`a $l$. Nous allons voir que
$\hlm$ et $\vlm$ sont des partitions.
Si $p=|\lambda|$ et $q$ est le nombre de cases qui ont un
nombre de crochets strictement sup\'erieur \`a $l$, alors 
$H^{p,q}[G(r,e),\ooo(l)]$ contient $S_\mu \C^e$, si $\mu$ est la partition
obtenue en r\'eordonnant toutes les parts $(l-\hlm_i)$ et $\vlm_j$.

Le seul inconv\'enient de cette m\'ethode est que l'ensemble des partitions
$l$-admissibles n'est pas ais\'e \`a d\'ecrire. Je propose dans ce paragraphe de
montrer qu'il est en bijection avec l'ensemble
de toutes les partitions de taille $(l-1,r)$, et d'\'etudier des propri\'et\'es
de cette bijection.

\begin{prop}
L'application $\lambda \mapsto \hlm$ est une bijection de
l'ensemble des partitions $l$-admissibles dans l'ensemble des partitions de
taille $(l-1,r)$.
\label{partition_admissible}
\end{prop}
\noindent
\dem
Tout d'abord si $\lambda$ est $l$-admissible, alors $\forall i,\hlm_i<l$
puisque les nombres de crochet sont strictement d\'ecroissants sur une ligne.
Introduisons quelques notations.
\begin{defi}
\ \\
\vspace{-.5cm}
\begin{itemize}
\item
Soit $h_{i,j}$ le nombre de crochet de la case situ\'ee sur la $i$-i\`eme
ligne et la $j$-\`eme colonne. Convenons que $h_{i,0}=+\infty$. 
\item
Pout tout $i$,
soit ${(g_{i,j})}_{j\in \N}$ la suite strictement croissante d'entiers dont 
l'image est le compl\'ementaire dans $\N$ de l'ensemble des $h_{i,j}$.
\item
Notons $\delta_i$ l'entier
$l-\hlm_i$, tel que $g_{i,\delta_i}=l$. Notons enfin $x_{i,j}$ le num\'ero de
colonne de la derni\`ere case sur la $i$-i\`eme ligne 
de nombre de crochet sup\'erieur \`a $g_{i,j}$.
\end{itemize}
\end{defi}
Ainsi, on a $g_{i,0}=0$, $x_{i,0}=\lambda_i$ et pour $j>>1$, $x_{i,j}=0$, par
nos conventions.

Alors, par la d\'efinition des nombres de
crochet, on a $h_{i-1,j}=h_{i,j}+\lambda_{i-1}-\lambda_i+1$. Les nombres
de crochet qui apparaissent sur la $(i-1)$-i\`eme ligne sont donc
$(1,\cdots, \lambda_{i-1}-\lambda_i$ et les 
$h_{i,j}+\lambda_{i-1}-\lambda_i+1$. Ainsi,
$g_{i-1,j+1}=g_{i,j}+\lambda_{i-1}-\lambda_i+1$. En
particulier, $g_{i-1,\delta_i+1}>l$, donc $\delta_{i-1}\leq \delta_i$ et donc
$\hlm$ est une partition. On d\'emontre maintenant un r\'esultat plus pr\'ecis
que la proposition \ref{partition_admissible}:

\begin{lemm}
Pout toute partition $\nu$ de taille $(l-1,r)$, il existe une unique
partition $l$-admissible $\lambda$ telle que $\hlm = \nu$; pour celle-ci, on a:
\begin{itemize}
\item{$\lambda_i=\lambda_{i+l-\nu_i} + \nu_i$.}
\item{$x_{i,j}=\lambda_{i+j}$.}
\end{itemize}
\label{inverse}
\end{lemm}
\dem
On d\'emontre par r\'ecurrence descendante sur $i_0$ qu'une partition a ses nombres
de crochet diff\'erents de $l$ \`a partir de la ligne $i_0$ et v\'erifie
$\hlm_i=\nu_i$ pour $i \geq i_0$ \ssi elle v\'erifie les deux points du lemme
pour $i\geq i_0$.

On a $\lambda_r=\nu_r,x_{r,0}=\lambda_r$ et $x_{r,j}=0$ pour $j>0$, de sorte
que
la propri\'et\'e de r\'ecurrence est vraie pour $i_0=r$. Soit maintenant $i_0$ fix\'e
tel que cette hypoth\`ese est v\'erifi\'ee pour $i_0+1$. Alors posons 
$\delta=\delta_{i_0}(=l-\nu_{i_0})$ et
$d=\lambda_{i_0}-\lambda_{i_0+1}$ ($d$ d\'epend de $\lambda_{i_0}$ qui est pour
l'instant inconnu). On a vu que $g_{i_0,0}=0$ et 
$g_{i_0,j+1}=g_{i_0,j}+d+1$ pour $j\geq 0$. Les nombres de 
crochet de $\lambda$ en-dessous et sur la $i_0$-i\`eme ligne sont diff\'erents de
$l$ \ssi $g_{i_0,\delta}=l$, ce qui \'equivaut donc \`a
$d=l-1-g_{i_0+1,\delta_{i_0-1}}$. 
Il existe donc exactement une possibilit\'e
pour $\lambda_{i_0}$. De plus, on a toujours
$x_{i_0,j+1}=x_{i_0+1,j}=\lambda_{i_0+j+1}$ et
le fait qu'effectivement $\hlm_{i_0}=\nu_{i_0}$
implique que
$\lambda_{i_0}=x_{i_0,\delta}+\nu_{i_0}=\lambda_{i_0+\delta}+\nu_{i_0}$.
\fin

\begin{nota}
Si $\nu$ est une partition de taille $(l-1,r)$, je noterai $\widehat \nu$
la partition $\lambda$ telle que $\hlm=\nu$. Par aileurs je noterai
$p(\nu)$ et $q(\nu)$ le poids et le nombre de cases de nombre de crochet
sup\'erieur \`a $l$ de $\lambda$.
\end{nota}

J'ai affirm\'e qu'aussi $\vlm$ est une partition, cela
d\'ecoule du fait que $\lambda^*$ est 
$l$-admissible et
$\vlm={}^h{(\lambda^*)}^-$.
\para

Si $\lambda_1,\lambda_2$ sont des partitions,
on note $\lambda_1 \subset \lambda_2$ le fait que 
$\forall i,\lambda_{1,i} \leq \lambda_{2,i}$.
Par une r\'ecurrence imm\'ediate, le lemme \ref{inverse} implique:
\begin{prop}
Si $\nu_1 \subset \nu_2$, alors 
$\widehat \nu_1 \subset \widehat \nu_2$.
\label{inclusion_inverse}
\end{prop}

\begin{prop}
Soient $k$ et $l$ deux entiers et $\lambda$ une partition. Si $\lambda$ est
$l$-admissible, alors elle est $(kl)$-admissible.
\end{prop}
\dem
En effet, si $G_i$ d\'esigne l'ensemble des entiers naturels qui
n'apparaissent pas parmi les nombres de crochet sur la $i$-i\`eme ligne, on sait
que pour tout $i$, $l\in G_i$ et il existe $d_i$ tel que 
$G_{i-1}=(d_i+1+G_i) \cup \{ 0 \}$.
Une r\'ecurrence imm\'ediate prouve alors que pour tout $i$, $G_i$ v\'erifie la
propri\'et\'e $x\in G_i \Rightarrow x+l\in G_i$. La proposition en d\'ecoule.
\fin

La contrapos\'ee de cette proposition implique que 
pour toute partition $\lambda$ 
exactement inscrite dans un rectangle
$r \times (e-r)$ (c'est-\`a-dire que $\lambda_1=e-r$ et $l(\lambda)=r$), 
il existe une case dont le nombre de crochet est $d$ pour tout
$d$ diviseur de $e-1$. Une telle partition n'est donc pas $d$-admissible. Si
$\lambda$ correspond \`a une composante non nulle de
$H^{*,*}[G(r,\C^e),\ooo(d)]$, on a donc soit $\lambda_1<e-r$ soit 
$l(\lambda)<r$. Ainsi, si $S_\mu\C^e \subset H^{*,*}[G(r,\C^e),\ooo(d)]$, on a
soit $\mu_e=0$, soit $\mu_1=d$.

%******************************************************************************

\subsection{Partitions admissibles de poids maximal}

Comme nous le verrons au paragraphe suivant, 
pour d\'emontrer des th\'eor\`emes d'annulation
pour la cohomologie de fibr\'es vectoriels amples, il est utile de bien
comprendre la cohomologie des fibr\'es homog\`enes sur une grassmannienne, et plus
pr\'ecis\'ement quelle est, pour $l,n$ et $r$ des entiers quelconques, la valeur
maximale de $p$ telle qu'il existe un entier $q$ avec
$H^{p,q}(G(r,n),\ooo(l))$ non nul. Si $n>rl$, cet entier a \'et\'e d\'etermin\'e
par Laurent Manivel \cite[proposition 1.2.1, p.111]{manivel_these}
et cela a permis de d\'emontrer de nombreux th\'eor\`emes
d'annulation, dont celui de F. Laytimi et W. Nahm concernant les 
partitions en forme de crochets \cite{nahm}. Lorsque $n,r$ et $l$ ne
satisfont plus cette in\'egalit\'e, la combinatoire des partitions $l$-admissibles
devient nettement plus compliqu\'ee et aucun r\'esultat ne donne cet entier $p$.

Soit $n,r,l$ des entiers fix\'es.
Si $a,\alpha,\beta,c,\gamma$ sont des entiers, notons 
$\mu(a,\alpha,\beta,c,\gamma)$ la partition $\mu$ d\'efinie par:
$$
\left \{
\begin{array}{ll}
\forall \ i\leq \alpha (l-a), & \mu_i=a\\
\forall \ \alpha (l-a)<i\leq \alpha (l-a) + \beta (l-a+1), & \mu_i=a-1\\
\forall \ \alpha (l-a) + \beta (l-a+1)<i\leq \alpha(l-a)+\beta(l-a+1)+\gamma, 
& \mu_i=c
\end{array}
\right .
$$

Cette partition est de longueur $(\alpha + \beta)(l-a) + \beta + \gamma$ et
un calcul direct utilisant le lemme
\ref{inverse} montre que 
$\widehat \mu_1(a,\alpha,\beta,c,\gamma)=(\alpha + \beta)a-\beta +c$.

Notons $\delta(\alpha,\beta)=\min\{\alpha,\beta\}$ si $\beta \not = 0$,
et $\delta(\alpha,0)=\alpha$. Rappelons que si $\mu$ est une partition,
$\widehat \mu$ a \'et\'e d\'efinie
comme l'unique partition telle que ${}^h\widehat \mu^-=\mu$ (cf lemme
\ref{inverse}).

\begin{theo}
Soient $n,r$ et $l$ des entiers, et $\lambda$ une partition $l$-admissible de
longueur $r$ et telle que $\lambda_1 \leq n-r$.\\
Alors, soit $\lambda$ est de la forme $\widehat \mu(a,\alpha,\beta,c,\gamma)$
avec $\gamma \leq l-a$,
soit il existe des entiers $a,\alpha,\beta,c,\gamma$ avec $\gamma \leq l-a$
et tels que
$\widehat \mu(a,\alpha,\beta,c,\gamma)$ soit aussi $l$-admissible de
longueur $r$ et de premi\`ere part inf\'erieure ou \'egale \`a
$n-r$, et que de plus
$|\widehat \mu(a,\alpha,\beta,c,\gamma)| \geq |\lambda| + 
\delta(\alpha,\beta)$. De plus, ces entiers sont tels que $\alpha + \beta$ et
$\gamma +c$ sont le quotient et le reste de la division euclidienne de $n$
par $l$.
\label{partition_p_max}
\end{theo}
\dem
Consid\'erons $n>r$ et $l$ trois entiers, et soit $\lambda$ une partition
$l$-admissible de longueur au plus $r$ et telle que $\lambda_1 \leq n-r$. 
On va raisonner
sur la partition $\hlm$: nous allons voir que dans trois cas, on peut modifier
cette partition pour en obtenir une de poids sup\'erieur et encore de longueur
au plus $r$ et de même premi\`ere part. D\'efinissons une suite $a_i$ par
$$
\left \{
\begin{array}{l}
a_1=1\\
a_{i+1}=\min \{ a_i+l-\hlm_{a_i} , r+1 \}.
\end{array}
\right .
$$
Soit $A$ le plus grand entier $i$ tel que $a_i \leq r$.
Le lemme \ref{inverse} montre que pour $i \leq A$,
$$\lambda_{a_i}=\sum_{i \leq j \leq A} \hlm_{a_j}.$$
En particulier, si $\mu$ est une partition telle que 
$\forall i\leq A,\mu_{a_i}=\hlm_{a_i}$, alors on a $\widehat \mu_1=\lambda_1$.
Par exemple, si $\mu$ la partition d\'efinie par
$\forall i\leq A, \forall j\in [a_i,a_{i+1}-1]$,
$\mu_j=\hlm_{a_i}$, on a, par la proposition \ref{inclusion_inverse},
$\widehat \mu \supset \lambda$, et
$\widehat \mu_1=\lambda_1$. Par exemple, pour $l=5$, on remplace la partition
$\yng(3,2,2,1,1,1)$ par $\yng(3,3,2,2,2,1)$.

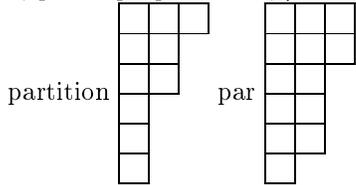
\lpara

Supposons qu'il existe des entiers $i,a$ et $b$ avec 
$b \leq a-2$ et tels que

\ \vspace{-3.3cm}

$$
\hspace{-3.5cm}
\begin{array}{cc}
\left \{
\begin{array}{ll} 
\forall j\in [i+1,i+l-a], & \hlm_j=a\\ 
\forall j\in [i+l-a+1,i+2(l-a)+1], & \hlm_j=b.
\end{array}
\right .

&

\setlength{\unitlength}{3mm}
\begin{picture}(5,12)(0,3)
\put(0,6){\line(1,0){8}}
\put(0,3){\line(1,0){8}}
\put(0,0){\line(1,0){4}}

\put(0,0){\line(0,1){6}}
\put(4,0){\line(0,1){3}}
\put(8,3){\line(0,1){3}}

\put(1.5,4.1){$a$}
\put(8.5,4.1){$l-a$}

\put(1.5,1.1){$b$}
\put(4.5,1.1){$l-a$}

\end{picture}

\end{array}
$$
\vspace{.5cm}

\noindent
Alors, consid\'erons la
partition $\mu$ d\'efinie par \vspace{-1cm}
 
$$
\begin{array}{cc}
\hspace{-2cm}
\left \{
\begin{array}{ll}
\forall j\in [i+1,i+l-a], & \mu_j=a-1\\
\forall j\in [i+l-a+1,i+2(l-a)+1], & \mu_j=b+1\\
\forall j \not \in [i+1,i+2(l-a)+1], & \mu_j=\hlm_j.
\end{array}
\right .

&

\setlength{\unitlength}{3mm}
\begin{picture}(8,6)(0,3)
\put(0,6){\line(1,0){8}}
\put(0,3){\line(1,0){8}}
\put(0,0){\line(1,0){4}}

\put(0,0){\line(0,1){6}}
\put(4,0){\line(0,1){3}}
\put(8,3){\line(0,1){3}}

\put(1,4.1){$a-1$}
\put(8.5,4.1){$l-a$}

\put(.5,1.1){$b+1$}
\put(4.5,1){$l-a$}

\end{picture}
\end{array}
$$

Le lemme \ref{inverse} implique que
$$
\left \{
\begin{array}{ll}
\forall j\in [i+1,i+l-a], & \widehat \mu_j=\lambda_j\\
\forall j\in [i+l-a+1,i+2(l-a)+1], & \widehat \mu_j\geq \lambda_j+1\\
\forall j \not \in [i+1,i+2(l-a)+1], & \widehat \mu_j=\lambda_j.
\end{array}
\right .
$$
Par cons\'equent, on a $|\widehat \mu| \geq |\lambda| + l-a + 1$. Appelons
``transformation A'' le fait de remplacer $\hlm$ par $\mu$.
\para

Supposons enfin l'existence d'entiers $i,a,b,\beta\geq 1$ tels que
\vspace{-1.5cm}

$$
\begin{array}{cc}
\hspace{-.8cm}
\left \{
\begin{array}{ll} 
\forall j\in [i+1,i+ l-a -1], & \hlm_j=a+1\\ 
\forall j\in [i+l-a,i+(1+\beta)(l-a)-1], & \hlm_j=a\\
\forall j\in [i+(1+\beta)(l-a),i+(2+\beta)(l-a)-1], & \hlm_j=b.
\end{array}
\right .

&

\setlength{\unitlength}{3mm}
\begin{picture}(8,6)(0,3)
\put(0,6){\line(1,0){4.5}}
\put(0,4){\line(1,0){4.5}}
\put(0,2){\line(1,0){4}}
\put(0,0){\line(1,0){1}}

\put(0,0){\line(0,1){6}}
\put(1,0){\line(0,1){2}}
\put(4,2){\line(0,1){2}}
\put(4.5,4){\line(0,1){2}}

\put(.7,4.6){$a+1$}
\put(5.1,4.6){$l-a-1$}

\put(1.5,2.6){$a$}
\put(5,2.6){$\beta(l-a)$}

\put(.2,.6){$b$}
\put(1.5,.6){$l-a$}

\end{picture}
\end{array}
$$
\vspace{.3cm}

\noindent
Alors, consid\'erons la
partition $\mu$ d\'efinie par \vspace{-1cm}
$$
\begin{array}{cc}
\hspace{-.8cm}
\left \{
\begin{array}{ll}
\forall j\in [i+1,i+ (1+\beta)(l-a) -1], & \mu_j=a\\
\forall j\in [i+(1+\beta)(l-a),i+(2+\beta)(l-a)-1], & \mu_j=b+1\\
\forall j \not \in [i+1,i+(2+\beta)(l-a)-1], & \mu_j=\hlm_j.
\end{array}
\right .

&

\setlength{\unitlength}{3mm}
\begin{picture}(8,6)(0,3)
\put(0,6){\line(1,0){4}}
\put(0,4){\line(1,0){4}}
\put(0,2){\line(1,0){4}}
\put(0,0){\line(1,0){1.5}}

\put(0,0){\line(0,1){6}}
\put(1.5,0){\line(0,1){2}}
\put(4,2){\line(0,1){4}}

\put(1.5,4.7){$a$}
\put(4.5,4.6){$l-a-1$}

\put(1.5,2.7){$a$}
\put(4.4,2.6){$\beta(l-a)$}

\put(.25,.8){\tiny{$b$}}
\put(.5,.8){\tiny{$+$}}
\put(.95,.8){\tiny{$1$}}
\put(2,.6){$l-a$}

\end{picture}
\end{array}
$$

Le lemme \ref{inverse} implique que
$$
\left \{
\begin{array}{ll}
\forall j\in [i,i+ l-a -1], & \widehat \mu_j=\lambda_j\\
\forall j\in [i+l-a,i+(2+\beta)(l-a)-2], & \widehat \mu_j = \lambda_j+1\\
\forall j \not \in [i+1,i+(2+\beta)(l-a)-1], & \widehat \mu_j=\lambda_j.
\end{array}
\right .
$$
En particulier, on a donc $|\widehat \mu| \geq |\lambda| + 2 (l-a)$. Appelons
``transformation B'' le fait de remplacer $\hlm$ par $\mu$.
\lpara

On peut donc construire une suite finie
de partitions $\mu^i$ par r\'ecurrence: tout d'abord, on pose
$\mu^0=\hlm$. Ensuite, on impose que pour passer
de $\mu^i$ \`a $\mu^{i+1}$, soit
on ajoute une case, soit on fait l'une des deux transformations A ou B, 
lorsque cela est possible de fa\c{c}on \`a ce que $\widehat \mu^{i+1}$ reste dans un
rectangle de taille $r \times (n-r)$. A un certain indice $N$, on ne peut
plus faire aucune de ces trois op\'erations. Soit $(b_i)$ les nombres d\'efinis
comme $a_i$ pour la partition $\mu^N$:
$$
\left \{
\begin{array}{l}
b_1=1\\
b_{i+1}=\inf \{ b_i+l-\mu^N_{b_i} , r+1 \}.
\end{array}
\right .
$$
Comme on ne peut plus ajouter de cases, on a
$\mu^N_j=\mu^N_{b_i}$ pour $b_i \leq j \leq b_{i+1}-1$, et comme aucune des
deux configurations \'etudi\'ees ne peut avoir lieu, $\mu^N$ est de type
$\mu(a,\alpha,\beta,c,\gamma)$ avec $\gamma \leq l-a$. De plus, puisque l'on ne
peut pas ajouter une case, on a soit $l[\mu(a,\alpha,\beta,c,\gamma)]=r$,
soit $c=0$. Si $c=0$, $\mu(a,\alpha,\beta,c,\gamma)$ ne d\'epend pas de $\gamma$;
on peut donc dans tous les cas supposer que
$$\alpha(l-a) + \beta(l-a+1) + \gamma = r.$$

Pour montrer que $|\mu^N| \geq |\lambda| + 
\delta(\alpha,\beta)$, il suffit de montrer que
$|\mu^N| \geq |\mu^{N-1}| + 
\delta(\alpha,\beta)$. Or, quand on applique l'op\'eration A, on obtient une
partition qui  a exactement $l-a$ parts \'egales \`a $a-1$, et quand on applique
l'op\'eration B, on obtient une partition qui a exactement $2(l-a)-1$ parts
\'egales \`a $a$: ce ne sont donc pas des partitions de type
$\mu(a,\alpha,\beta,c,\gamma)$. On en d\'eduit donc que pour passer de
$\mu^{N-1}$ \`a $\mu^N=\mu(a,\alpha,\beta,c,\gamma)$, on a ajout\'e une case. Soit
cette case se trouvait sur la $r$-i\`eme ligne, et on a
$|\mu^N| - |\mu^{N-1}| = \alpha + \beta + 1$; soit, si $\beta>0$, elle se
trouvait sur la $\alpha(l-a)+\beta(l-a+1)$-i\`eme ligne, et on a
$|\mu^N| - |\mu^{N-1}| = \beta$; soit, enfin, elle se trouvait sur la
$\alpha(l-a)$-i\`eme ligne, et dans ce cas
$|\mu^N| - |\mu^{N-1}| = \alpha$; on voit que dans tous les cas,
$|\mu^N| - |\mu^{N-1}| \geq \delta(\alpha,\beta)$.

Pour achever la preuve du th\'eor\`eme, il ne reste plus qu'\`a calculer 
$\alpha + \beta$ et $\gamma + c$. Or, on a vu que
$\widehat \mu_1(a,\alpha,\beta,c,\gamma)=(\alpha + \beta)a-\beta +c$;
comme les transformations effectu\'ees ne changent pas la premi\`ere part, on en
d\'eduit $\lambda_1=(\alpha+\beta)a - \beta + c$. Par ailleurs,
on a $r=(\alpha + \beta)(l-a)+ \beta + \gamma$. On en d\'eduit donc que
\begin{equation}
n=r+\lambda_1=(\alpha + \beta)\ l+(\gamma +c).
\label{division_euclidienne}
\end{equation}
Or, on a vu que $\gamma \leq l-a$ et $c \leq a$. Par convention, on exclut
les \'egalit\'es $\gamma = l-a$ et $c=a$,
car $\mu(a,\alpha,0,a,l-a)=\mu(a,\alpha+1,0,0,0)$; on a donc
$\gamma + c < l$. La formule (\ref{division_euclidienne}) exprime donc la
division euclidienne de $n$ par $l$, concluant la preuve du th\'eor\`eme.
\fin

%******************************************************************************

%******************************************************************************

\section{Suite spectrale de Borel Le-Potier}

Dans ce paragraphe, je rappelle la d\'efinition de la suite spectrale de
Borel Le-Potier, qui est \`a la base de mes th\'eor\`emes d'annulation.

Soit $\pi:Y \rightarrow X$ une fibration localement triviale
propre et $E$ un fibr\'e vectoriel sur $Y$.
Soit $\Omega^i_{Y/X}$ le fibr\'e des $i$-formes sur $Y$ relatives \`a $\pi$, 
d\'efini par $\Omega^i_{Y/X} = \Lambda^i \Omega_{Y/X}$ et
la suite exacte de fibr\'es sur $Y$ suivante:
$$ 
0 \rightarrow \pi^*\Omega X \stackrel{\pi^*}{\rightarrow} \Omega Y \rightarrow 
\Omega_{Y/X} \rightarrow 0.
$$
Soient aussi les fibr\'es
$G^{t,p}:=\Omega^{p-t}_{Y/X} \otimes \pi^* \Omega^t_X$.
Pour chaque $p$, on a \cite{lepotier1} une suite spectrale 
aboutissant sur $H^{p,q}(Y,E)$ et dont les
termes d'ordre 1 sont:
$$ \Es{p}{t}{q-t}{1} = H^q (Y,G^{t,p} \otimes E). $$

\lpara

\label{para_leray}
Pour calculer les groupes de cohomologie $H^q (Y,G^{t,p} \otimes E)$, on
utilise une suite spectrale de Leray. La suite spectrale de termes d'ordre 2
\begin{equation}
\Ed{p}{t}{k}{j-k}{2} = H^k (X,R^{j-k}\pi_* G^{t,p}) 
= H^{t,k}[X,R^{j-k}\pi_* (E \otimes \Omega^{p-t}_{Y/X})]
\label{leray}
\end{equation}
aboutit sur $\Es{p}{t}{j}{1}$. Introduisons enfin une notation:
\begin{nota}
Si $\pi:Y \rightarrow X$ est une fibration et $E$ un fibr\'e vectoriel sur $Y$,
notons $R^{p,q} \pi_* E$ le faisceau $R^q \pi_* (E \otimes \Omega^p_{Y/X})$.
\label{rpq}
\end{nota}

%******************************************************************************

%******************************************************************************

\section{Th\'eor\`eme d'annulation pour un produit tensoriel de crochets}

W. Nahm et F. Laytimi ont montr\'e un th\'eor\`eme d'annulation pour la cohomologie
d'une puissance de Schur correspondant \`a un crochet d'un fibr\'e ample
\cite[th 2.1]{nahm}. Dans ce
paragraphe, je g\'en\'eralise leur r\'esultat \`a un produit de crochets. Si $\alpha$
et $\beta$ sont des entiers, notons $Z^{\alpha,\beta}$ le foncteur
$S_\lambda$ pour la partition $\lambda$ de longueur $\alpha+1$ et de poids 
$\beta$ telle que $\lambda_1=\beta-\alpha$ et 
$\forall 1<i\leq \alpha + 1, \lambda_i=1$.

\begin{prop}
Soit $E$ un fibr\'e ample de rang $e$ sur une vari\'et\'e $X$
projective et lisse de 
dimension $n$. Soit $a\in \N$,
$(k_i)_{1 \leq i \leq a}$ et $(\alpha_i)_{1 \leq i \leq a}$ des entiers tels
que $\alpha_i < k_i$; soit $\sigma=\sum \alpha_i$ et $k=\sum k_i$; alors,
$$ H^{p,q} (X, \otimes_i Z^{k_i-\alpha_i-1,k_i} E)=0 \mbox{ si } 
q> n-p+[\delta(n-p)+\sigma][ae-k+2\sigma] - \sigma(\sigma+1). $$
\label{proposition_q}
\end{prop}
\vspace{-.7cm}
Posons $Q(p,\sigma)=n-p+[\delta(n-p)+\sigma][ae-k+2\sigma] - \sigma(\sigma+1)$.
\\
\dem
Reprenant les id\'ees de \cite{nahm},
soit $C_2^k=\frac{k(k-1)}{2}$ et $\delta$ la fonction d\'efinie par
$\forall n \in \N, C^{\delta(n)}_2 \leq n < C^{\delta(n)+1}_2$.
On peut alors d\'efinir un ordre sur $\N^2$.
Celui-ci est donn\'e par
l'ordre lexicographique sur $\N^3$ et l'injection 
$(x,\sigma) \mapsto (\delta(x)+\sigma,x-C^{\delta(x)}_2,\sigma)$. 
Notons $\NN$ l'ensemble $\N^2$ muni de cet ordre. L'int\'er\^et
principal de $\NN$ est qu'il v\'erifie:

\begin{lemm}\cite{nahm}:
Pour tous $(x,\sigma)\in \NN,\mu \in \Z-\{ 0 \}$, si $r=\delta(x)$
et si $x+\mu r \in \N$, alors
$(x+\mu r,\sigma-\mu)<(x,\sigma)$. \label{ordre}
\end{lemm}

On d\'emontre la proposition par r\'ecurrence sur les couples
$(n-p,\sigma) \in \NN$.
Cette r\'ecurrence peut sembler
peu naturelle et troubler le lecteur; aussi, je donne les 10 premiers \'el\'ements
de $\NN$ pour qu'il puisse suivre plus ais\'ement les premiers pas
de la r\'ecurrence.

$$
\begin{array}{cccccc}
 & \ n-p\ & \ \sigma \ & \ \delta(n-p) + \sigma\ 
& \ n-p-C_2^{\delta(n-p)}\ & \ \sigma \ \\
1 & 0 & 0 & 1 & 0 & 0 \\
2 & 1 & 0 & 2 & 0 & 0 \\
3 & 0 & 1 & 2 & 0 & 1 \\
4 & 2 & 0 & 2 & 1 & 0 \\
5 & 3 & 0 & 3 & 0 & 0 \\
6 & 1 & 1 & 3 & 0 & 1 \\
7 & 0 & 2 & 3 & 0 & 2 \\
8 & 4 & 0 & 3 & 1 & 0 \\
9 & 2 & 1 & 3 & 1 & 1 \\
10 & 5 & 0 & 3 & 2 & 0 \\
\end{array}
$$

Si on place ces num\'eros d'apparition sur un plan rep\'er\'e par $n-p$ selon l'axe 
des abscisses et $\sigma$ selon l'axe des ordonn\'ees, on obtient le sch\'ema
suivant:
$$
\begin{array}{ccccccc}
 & & & & & & \sigma \\
 & & & & & & \uparrow \\
 & & & & & & 7\\
 & & & & 9 & 6 & 3\\
n-p \ \ \leftarrow & 10 & 8 & 5 & 4 & 2 & 1\\
\end{array}
$$

\para

La m\'ethode employ\'ee est celle maintenant classique qui consiste \`a constater que
les groupes de cohomologie que l'on cherche \`a annuler apparaissent dans une
suite spectrale qui calcule la cohomologie d'un fibr\'e en droites ample sur une
vari\'et\'e ad\'equate $\cal Y$ elle-m\^eme fibr\'ee au-dessus de $X$. 
Ce dernier groupe est nul par
le th\'eor\`eme de Kodaira si $q$ est suffisant; il suffit donc de s'assurer que
ces groupes ``passent \`a travers'' 
la suite spectrale, ce qui fait intervenir une
r\'ecurrence. 

Soit donc $X^n,E^e,p,q,(k_i)$  et $\sigma_0$ v\'erifiant les hypoth\`eses de la
proposition fix\'es. On va montrer l'annulation
de $H^{p,q}(X,\produit{\sigma_0}{E})$ en supposant 
la proposition vraie pour tous
les couples $(p',\sigma)$ tels que $(n-p',\sigma) < (n-p,\sigma_0)$. Supposons 
la suite 
$(k_i)$ croissante. Soit $r=\delta(n-p)$; soit aussi
$l_i$ et $s_i$ les quotients et restes de la division de $k_i$ par $r$:
$k_i=rl_i+s_i$. Un premier
lemme assure que l'on peut supposer que $e>r$ tous les $l_i$ sont
sup\'erieurs ou \'egaux \`a 
$$l_0=
\left \{
\begin{array}{lll}
\frac{r\sigma_0 + (n-p)}{r-1} & si & r>2 \\
\sigma + 1                & si & r=1
\end{array}
\right . $$

\begin{lemm}
Soit $l$ un entier tel que la proposition \ref{proposition_q}
soit vraie pour une 
certaine valeur de $n-p$ et
de $\sigma$, et pour toutes les vari\'et\'es et tous les fibr\'es vectoriels amples,
\`a la condition que tous les $l_i$ soient sup\'erieurs ou \'egaux \`a $l$ et que
$e>\delta(n-p)$. Alors
cette proposition est vraie sans restrictions pour ces valeurs de $e$ et de 
$\sigma$.
\label{astuce_nahm}
\end{lemm}
\dem
Soit $E$ un fibr\'e sur une vari\'et\'e $X$ et $L$ un fibr\'e en droites ample
sur $X$. On peut supposer $l>\delta(n-p)$. 
La proposition est vraie pour $l_i'=l_i + l$
(c'est-\`a-dire, puisque $k_i'=rl_i'+s_i$, pour $k_i'=k_i+rl$)
et les fibr\'es amples $E'$ de rang
sup\'erieur \`a $\delta(n-p)$. Pour le 
fibr\'e ample 
$E'=E\oplus L^{\oplus rl}$ de rang $e+rl$, et pour $k_i'=k_i+rl$,
on a donc l'annulation de 
$H^{p,q} (X, \otimes_i Z^{k_i'-\alpha_i-1,k_i'} E')=0$ si 
$$
\begin{array}{rcl}
q & > & [\delta(n-p)+\sigma][a(e+rl)-k'+2\sigma] - \sigma(\sigma+1)\\
  & = & [\delta(n-p)+\sigma][ae-k+2\sigma] - \sigma(\sigma+1).
\end{array}
$$ 
Mais comme
$Z^{k_i'-\alpha_i-1,k_i'} E' \supset Z^{k_i-\alpha_i-1,k_i} E \otimes 
L^{\otimes rl}$, on en d\'eduit que pour toute vari\'et\'e $X$, tous fibr\'es $E$ et
$L$, $H^{p,q} (X, \otimes_i Z^{k_i-\alpha_i-1,k_i} E \otimes 
L^{\otimes arl})=0$ si 
$q> [\delta(n-p)+\sigma][ae-k+2\sigma] - \sigma(\sigma+1)$.
En vertu du lemme suivant, ce r\'esultat reste
vrai si $L$ est trivial, et le lemme \ref{astuce_nahm} est donc prouv\'e.
\fin

\begin{lemm}
Soient $n,p,q_0,e$ des entiers et $\lambda$ une partition tels que
$H^{p,q}(X,S_\lambda E \otimes L)=0$ pour toute vari\'et\'e projective lisse $X$
de dimension $n$, 
tout fibr\'e vectoriel ample $E$ de rang $e$
et tout fibr\'e en droites ample $L$ sur $X$, si $q>q_0$.

Alors $H^{p,q}(X,S_\lambda E \otimes L)=0$
sous les mêmes conditions, sauf que $L$ est seulement nef.
\end{lemm}
\dem
La d\'emonstration du lemme 1.5.1 p.128 dans \cite{manivel}
est valable sous ces hypoth\`eses.
\fin

On suppose donc dor\'enavant que $l_1 \geq l_0$ et $e>r$. Si
$0<s<r<e$ sont des entiers et $V$ un espace vctoriel de dimension $e$, 
on notera $M_{s,r}(V)$ la vari\'et\'e de drapeaux absolue constitu\'ee de l'ensemble
des couples $(V_r,V_s)$ de sous-espaces vectoriels de $V$ avec
$$ 0 \subset V_r \subset V_s \subset V, \ \ codim(V_r)=r, \ \ codim(V_s)=s. $$
Soit alors ${\cal M}_{s,r}(E)$ la vari\'et\'e de drapeaux relative \`a $E$, 
c'est-\`a-dire la vari\'et\'e 
fibr\'ee au-dessus de $X$ dont la fibre au-dessus du point $x$ s'identifie \`a la
vari\'et\'e $M_{s,r}(E_x)$ des drapeaux de la forme
$$ (0 \subset E_r \subset E_s \subset E_x), \ \ codim(V_r)=r, \ \ 
codim(V_s)=s. $$
Soit aussi ${\cal Q}^{l+1,l}$ le fibr\'e en droites sur ${\cal M}_{s,r}(V)$ 
dont la fibre au-dessus
du drapeau pr\'ec\'edent s'identifie \`a 
$(\det (E_x/E_s))^{l+1} \otimes (\det (E_s/E_r))^l = \det (E/E_s) \otimes
(\det (E/E_r))^l$. Consid\'erons alors le produit au-dessus de $X$ d\'efini par
${\cal Y} := \times M_{s_i,r} (E) $ et le fibr\'e en droites $\cal L$ au-dessus 
de $\cal Y$ \'egal au
produit $\otimes \pi_i^* {\cal Q}_i$, 
si $\pi_i$ d\'esigne la projection de $\cal Y$
sur ${\cal M}_{s_i,r_i} (E)$ et ${\cal Q}_i$ le fibr\'e ${\cal Q}^{l_i+1,l_i}$ 
sur cette vari\'et\'e de
drapeaux relative. Comme $e>r>s_i$ pour tout $i$, cette vari\'et\'e a bien un sens.

\para

Laurent Manivel a \'etudi\'e une partie de la cohomologie des fibr\'es en droites
$Q^{l+1,l}$ sur une vari\'et\'e de drapeaux absolue, 
partie suffisante pour \'etablir notre
proposition. N\'eanmoins, j'ai besoin de 
g\'en\'eraliser ses r\'esultats \`a un produit de
vari\'et\'es de drapeaux. Cette g\'en\'eralisation sera une cons\'equence facile de la
formule de Kunneth.

\para

Soit donc $M_{s,r}(V)$ une vari\'et\'e de drapeaux absolue et $Q^{l+1,l}$ comme 
pr\'ec\'edemment. Notons $\lambda=l-1,t(r)=\frac{r(r+1)}{2}$ et $k=rl+s$. 
Notons $\pi_{r,s}=\frac{s(2r-s+1)}{2}$ et pour $\pi$ un entier, 
$$n_s(\pi)=\# \{ c \in \{ 0,1 \}^r, |c|=s \mbox{ et } \sum ic_i = \pi \}$$
Ainsi $n_s(\pi)=0$ si $\pi<0$ ou $\pi>\pi_{r,s}$.

Alors, on a 
\cite[proposition 1.2.1, p.111 et lemme 1.3.1, p.114]{manivel_these}:
\begin{prop}
Supposons $p=\lambda . t(r) + \pi$ avec 
$\pi \geq \pi_{r,s} - k + l$. Alors
$$
H^{p,q}(M_{s,r}(V),Q^{l+1,l})= \bigoplus_{\alpha=0}^l 
n_s(\pi + r\alpha) \delta_{q,p-r\lambda-s+\alpha} Z^{k-\alpha-1,k}V. $$
\label{manivel_drapeaux}
\end{prop}

Consid\'erons maintenant un produit $Y=\times M_{s_i,r}(V)$ de vari\'et\'es de 
drapeaux et le fibr\'e en droites $L=\otimes \pi_i^* Q_i$, avec $Q_i$ le
fibr\'e $Q^{l_i+1,l_i}$ sur $M_{s_i,r}(V)$ et $\pi_i$ la projection naturelle
de $Y$ sur $M_{s_i,r}$. Je propose alors la g\'en\'eralisation suivante de la
proposition \ref{manivel_drapeaux}, si $\lambda_i$ d\'esigne $l_i-1$, 
$\lambda=\sum \lambda_i$, $l=\sum l_i$, $s=\sum s_i$, et $k=\sum k_i$: 

\begin{prop}
Supposons que $p=\lambda.t(r) + \pi$ avec 
$\pi \geq \sum \pi_{r,s_i} - k + l$.\\
Soit $\sigma=q+r\lambda+s-p$. Alors
$$ H^{p,q}(Y,L)= n_s(\pi + r \sigma)
\bigoplus_{\alpha_i:\sum \alpha_i=\sigma} \otimes_i Z^{k_i-\alpha_i-1,k_i}V. $$
\label{produit_drapeaux}
\end{prop}
\vspace{-.4cm}
\noindent
Dans cette proposition, $n_s(\pi)$ d\'esigne le nombre
$$n_s(\pi):=\sum_{\sum \pi_i = \pi} n_{s_i}(\pi_i)= \# \{c_{i,j} 
\in \{0,1\}^{r^2}:
\forall i \ |c_i|=s_i \mbox{ et } \sum_{i,j} j c_{i,j}=\pi \}.$$
Cette proposition montre que le fait qu'un groupe 
$\otimes_i Z^{k_i-\alpha_i-1,k_i}V$ soit ou non dans
$ H^{p,q}(Y,L)$ ne d\'epend que de la somme $\sigma$ des $\alpha_i$, et non
des valeurs de tous les $\alpha_i$; c'est ce qui fait que l'on obtient un
th\'eor\`eme d'annulation o\`u la borne ne d\'epend aussi que de $\sigma$.
\para
\dem 
On applique tout d'abord la formule de Kunneth:
$$ H^{p,q}(Y,L)=\bigoplus_{\sum p_i=p,\sum q_i=q} \otimes 
H^{p_i,q_i} (M_{s_i,r},Q^{l_i,l_i+1}). $$
Ecrivons chaque $p_i$ comme $\lambda_i t(r) + \pi_i$.
Si une telle suite d'entiers $p_i$ fournit un groupe non nul, alors par la
proposition \ref{manivel_drapeaux}, on a pour tout $i$ la relation
$\pi_i \leq \pi_{r,s_i}$. Si donc il existe $i$ tel que 
$\pi_i < \pi_{r,i} - k_i + l_i$, comme $\pi=\sum \pi_i$, on en d\'eduit que
$\pi < \sum \pi_{r,_i} - k + l$, ce qui contredit l'hypoth\`ese de la 
proposition. Ainsi, pour
tout $i$ on a $\pi_{r,s_i} - k_i + l_i \leq \pi_i \leq \pi_{r,s_i}$, et donc
on peut en d\'eduire  par la proposition \ref{manivel_drapeaux} que
$$H^{p_i,q_i} (M_{s_i,r},Q^{l_i,l_i+1})=\oplus_{\alpha_i} n_s(\pi_i + r 
\alpha_i ) \delta_{q_i,p_i-r\lambda_i-s_i+\alpha_i} Z^{k_i-\alpha_i-1,k_i}V.$$

Si donc $\otimes_i Z^{k_i-\alpha_i-1,k_i}V$ apparaît dans notre groupe de
cohomologie, cela implique
que $q_i=p_i-r\lambda_i-s_i+\alpha_i$, et par
sommation que $q=p-r\lambda-s+\sigma$. Enfin, la multiplicit\'e de ce groupe
est bien la somme sur les $\pi_i$ des produits des multiplicit\'es de
$Z^{k_i-\alpha_i-1,k_i}V$ dans $H^{p_i,q_i} (M_{s_i,r},Q^{l_i,l_i+1})$, soit
$n_s(\pi)$.
\fin

Notons $\PM=\lambda.t(r)+\sum \pi_{r,s_i}$ et 
$\QM=\lambda.t(r-1)+\sum \pi_{r,s_i}-S$. Une cons\'equence de la proposition
\ref{produit_drapeaux} est:

\begin{prop}
Soit $p$ et $q$ deux entiers. Alors:
\begin{itemize}
\item{Si $p>\PM$ ou $q>\QM$, alors $H^{p,q}(Y,L)=0$.}
\item{Si $p\geq \PM-r\sigma$ ou $q\geq \QM-(r-1)\sigma$, alors
$$H^{p,q}(Y,L) \subset \bigoplus_{\sum \alpha_i \leq \sigma} \otimes_i
Z^{k_i-\alpha_i-1,k_i} V.$$}
\end{itemize}
\label{majoration_pq}
\end{prop}

\Para

On consid\`ere maintenant la fibration en vari\'et\'es de drapeaux
$\pi:{\cal Y} \rightarrow X$
introduite pr\'ec\'edemment. Soit $\Es{P}{i}{j}{m}$
la suite spectrale de Borel-Le Potier, avec
$P=p+\PM-r\sigma_0$.
Le r\'esultat
souhait\'e va \^etre cons\'equence de propri\'et\'es de cette suite spectrale. La 
premi\`ere exhibe un terme de la suite spectrale de Borel-Le Potier qui contient
le groupe que l'on veut annuler:

\begin{lemm}
Soit $q_0=\QM - (r-1)\sigma_0$. Alors
$$\Es{P}{p}{q+q_0-p}{1} \supset H^{p,q}(X,\produit{\sigma_0}{E}).$$
\label{notre_groupe}
\end{lemm}
\dem 
Comme on a suppos\'e que 
$l_1(r-1) \geq r\sigma_0 + (n-p)$, on a
$$
\begin{array}{rcl}
n + \PM - k_1 + l_1     &   =   & n + \PM - l_1(r-1) - s_1\\
                        & \leq  & n + \PM - r \sigma_0 - (n-p)\\
                        &   =   & \PM + p - r\sigma_0 = P.
\end{array}
$$
Ainsi, si $p \leq n$, alors
$P-p \geq P-n \geq \PM - k_1 + l_1$. Si $p>n$, il est clair que
$\Es{P}{p}{q+q_0-p}{1}=0$. On peut donc utiliser les 
propositions \ref{produit_drapeaux} et \ref{majoration_pq} pour obtenir des
renseignements sur les
$H^{P-p,q'}(Y,L)$ pour tous $q'$.
Pour tous les entiers $q_1$ et $q_2$, on a (cf notation \ref{rpq})
$$\Ed{P}{p}{q_1}{q_2}{2}=H^{p,q_1}[X,R^{P-p,q_2}\pi_* {\cal L}].$$
Si $f_{i,j}$ sont les applications de changement de cartes de $E$, on a vu que
$R^{P-p,q_2}\pi_* {\cal L}$ est un fibr\'e vectoriel de fibre type
$H^{P-p,q_2} (Y,L)$ (soit $x\in X$; $Y=\pi^{-1}(x)$ 
et $L$ est la restriction du fibr\'e en droites ${\cal L}$ \`a $Y$), et dont les 
changements de carte sont induits par $f_{i,j}$. Supposons que
$H^{P-p,q_2} (Y,L)=S_\lambda E_x$. Nous savons que $|\lambda|=\sum rl_i+s_i$.
On en d\'eduit donc que pour $g \in GL(E_x)$, l'application induite par $g$
sur $H^{P-p,q_2} (Y,L)$ est $S_\lambda g$. En effet, c'est vrai $g \in SL(E_x)$
par le th\'eor\`eme de Bott et pour $g$ une homoth\'etie, puisque l'application
induite par $\lambda.Id$ est $\lambda^{\sum rl_i+s_i}$ (l'action de
$\lambda.Id$ sur le fibr\'e tangent est triviale).

Les applications de changement de cartes de $R^{P-p,q_2}\pi_* {\cal L}$ sont
donc les applications $S_\lambda f_{i,j}$, si $f_{i,j}$ est une application
de changement de cartes de $E$. On a donc que
$R^{P-p,q_2}\pi_* {\cal L}=S_\lambda E$.
\para

Or
$H^{P-p,q_0}(Y,L) = \produit {\sigma_0}{V}$. En effet, si un
terme de la forme\\ $\produit {\sigma}{V}$ 
est une composante de ce groupe, alors on doit avoir
$\sigma = q_0 + rL + S - (P-p) = \sigma_0$. Par ailleurs, il est clair que 
$n_s(\sum \pi_{r,s_i})=1$. On en d\'eduit donc que 
$\Ed{P}{p}{q}{q_0}{2}=H^{p,q}(X,\produit{\sigma_0}{E})$.

Par ailleurs, si $q_2>q_0$ et $q_1$ est quelconque, 
$\Ed{P}{p_1}{q_1}{q_2}{2}=0$,
toujours par la proposition \ref{majoration_pq}. 
Enfin, si $q_1>q$ et $q_2<q_0$,
et que 
$H^{p,q_1}(X,\produit{\sigma}{E}) \subset \Ed{P}{p_1}{q_1}{q_2}{2}$, 
on a $\sigma < \sigma_0$ et $H^{p,q_1}(X,\produit{\sigma}{E})=0$ par 
l'hypoth\`ese de r\'ecurrence. Toutes les diff\'erentielles $d_m$ issus de ou 
aboutissant sur $\Ed{P}{p}{q_0}{q}{m}$ sont donc nulles, et
$\Ed{P}{p}{q_0}{q}{\infty}=H^{p_1,q_1}(X,\produit{\sigma_0}{E})$, ce qui 
implique notre lemme.
\fin

Pour voir que notre groupe pr\'ec\'edent passe
\`a travers la suite spectrale de Borel-Le Potier, il suffit de  montrer 
l'annulation des groupes $\Es{P}{p+m}{q-p-m+\sgn (m)+q_0}{1}$, pour tous les
entiers $m$ non nuls. Soit donc $m$ un tel entier et $p_1,q_1$ et $q_2$ des
entiers tels que $\Ed{P}{p_1}{q_1}{q_2}{\infty}$ soit un \'el\'ement d'une
filtration de $\Es{P}{p+m}{q-p-m+\sgn (m)+q_0}{1}$, c'est-\`a-dire des entiers
tels que:\\
$p_1=p+m \mbox{ et }q_1+q_2-p_1=q-p-m+\sgn (m)+q_0$.
Remarquons que de mani\`ere \'equivalente, on a
$p_1=p+m$ et $q_1+q_2=q+q_0+\sgn (m)$.

De nouveau, il suffit de montrer que tous les groupes
$\Ed{P}{p_1}{q_1}{q_2}{2}$ sont nuls. Or ceux-ci valent
$H^{p_1,q_1}[X,H^{P-p_1,q_2}({\cal Y},{\cal L})]$. Supposons que
$$\produit{\sigma}{V} \subset H^{P-p_1,q_2}(Y,L).$$ Alors par la proposition
\ref{majoration_pq}, $P-p_1 \leq \PM - r\sigma =(P-p) - r(\sigma-\sigma_0)$,
donc $p_1 \geq p + r(\sigma-\sigma_0)$, soit $m \geq  r(\sigma-\sigma_0)$. Par
ailleurs cette m\^eme proposition assure que
$q_2 \leq \QM - (r-1) \sigma = q_0 - (r-1)(\sigma - \sigma_0)$.

L'\'egalit\'e $q_1 + q_2 = q + q_0 + \sgn (m)$ donne alors
$q_1 \geq q + \sgn(m) + (r-1)(\sigma - \sigma_0)$. Notons 
$\mu=\sigma - \sigma_0$, comme nous avons vu que $m\geq r\mu$, 
$\sgn (m) \geq \sgn (\mu)$ et nous avons donc \'etabli:
$$
\begin{array}{c}
\sigma = \sigma_0 + \mu\\
p_1 \geq p + \mu r\\
q_1 \geq q + \sgn (\mu) + (r-1)\mu.\\
\end{array}
$$
Sous ces hypoth\`eses, il ne nous reste plus, en utilisant l'hypoth\`ese de 
r\'ecurrence, qu'\`a prouver l'annulation de $H^{p_1,q_1}(X,\produit{\sigma}{E})$.
Remarquons tout d'abord que l'on peut supposer que $ae-k+\sigma \leq 0$, 
car sinon $\produit{\sigma}{E}=0$.

\Para
\begin{lemm}
Dans l'ensemble ordonn\'e $\NN$, on a:
$(n-p_1,\sigma) < (n-p,\sigma_0)$
\label{ok_recurrence_q}
\end{lemm}
\dem 
Cet ordre \'etant croissant selon les deux coordonn\'ees (c'est-\`a-dire que si
$x\leq x',y \in \N$, alors $(x,y) \leq (x',y)$ et $(y,x) \leq (y,x')$), cela
d\'ecoule du lemme \ref{ordre}.
\fin

\begin{lemm}
Si $q>Q(p,\sigma_0)$ alors $q_1>Q(p_1,\sigma)$.
\label{q_suffisant}
\end{lemm}
\dem 
On peut supposer que $p_1=p+\mu r$ et $q_1=q+\mu (r-1)+\sgn (\mu)$. 
Calculons alors
$Q(p,\sigma_0) - Q(p+\mu r,\sigma_0+\mu) +\sgn (\mu) + (r-1)\mu$. 
Cette quantit\'e vaut:\\
$2\mu r + 2 \mu \sigma_0 + \mu^2 + \sgn (\mu)+ (r+\sigma_0)(ae-k+2\sigma_0)
-[\delta(n-p-\mu r) + \sigma_0 + \mu][ae-k + 2\sigma_0 + 2\mu]$. En remarquant
que $2\mu r + 2 \mu \sigma_0 =2\mu(r+\sigma_0)$ et en factorisant par
$r + \sigma_0$, on en d\'eduit que cette quantit\'e \'egale\\
$(ae-k+2\sigma_0+2\mu)(r+\sigma_0)-(ae-k+2\sigma_0-2\mu)
[\delta(n-p+\mu r) + \sigma_0 - \mu] + \mu^2 + \sgn(\mu)$,
soit $(ae-k+2\sigma)[r+\mu-\delta(r+\mu r)]+\mu^2+\sgn (\mu)$.

Comme $ae-k+\sigma \geq 0$, il d\'ecoule de la d\'efinition de $\delta$
que ce nombre est positif.
\fin

Enfin, on a l'in\'egalit\'e suffisante pour assurer que 
$\Es{P}{p}{q-p+q_0}{\infty}=0$:
\begin{lemm}
On a $P+q+q_0-\dim {\cal Y}>0$.
\label{E_infini_nul}
\end{lemm}
\dem
Tout d'abord, 
$$
2\sum \pi_{r,s_i}-S = 2rS + S - \sum s_i^2 -S \geq rS \mbox{ et }k=(L+a)r+S.
$$
Ainsi:
$$
\begin{array}{cl}
  & P+q+q_0-\dim {\cal Y}\\ 
= & p + q + Lt(r) + Lt(r-1) + 2 \sum \pi_{r,s_i} - S
- (2r-1)\sigma_0 - n - ar(e-r) \\
 \geq & \sigma_0(ae-k+\sigma_0)
\end{array}
$$

\fin 
\Para

Nous avons ainsi (p\'eniblement) achev\'e la 
preuve de la proposition \ref{proposition_q}.
Nous allons maintenant montrer une autre proposition qui
r\'etablit la sym\'etrie entre $p$ et $q$:

\begin{prop}
Soit $E$ un fibr\'e ample de rang $e$ sur une vari\'et\'e $X$ 
projective lisse et de dimension $n$. Soit 
$a\in \N$,
$(k_i)_{1 \leq i \leq a}$ et $(\alpha_i)_{1 \leq i \leq a}$ des entiers tels
que $\alpha_i < k_i$; soit $\sigma=\sum \alpha_i$ et $k=\sum k_i$, alors
$$ H^{p,q} (X, \otimes_i Z^{k_i-\alpha_i-1,k_i} E)=0 \mbox{ si } 
p > n - q + [\delta(n-q)+\sigma][ae-k+2\sigma] - \sigma(\sigma+1). $$
\label{proposition_p}
\end{prop}
\vspace{-.6cm}
\dem
Posons $P(q,\sigma)=n-q+[\delta(n-q)+\sigma][ae-k+2\sigma] - \sigma(\sigma+1)$.
La d\'emonstration est tr\`es similaire \`a celle de la proposition
\ref{proposition_q}. En effet, soit $r=\delta(n-q)$ et comme pr\'ec\'edemment
${\cal Y} := \times M_{s_i,r} (E) $ et le fibr\'e en droites $\cal L$ au-dessus 
de $\cal Y$ \'egal au
produit $\otimes \pi_i^* {\cal Q}_i$ , o\`u $l_i$ et $s_i$ sont 
le quotient et le
reste de la division euclidienne de $k_i$ par $r$ et 
${\cal Q}_i={\cal Q}^{l_i,l_i+1}$. Soit aussi
$p+P_{max}-r\sigma_0$, le lemme \ref{notre_groupe}
reste vrai et pour voir que ce groupe passe \`a travers la suite spectrale, il
suffit comme pr\'ec\'edemment de montrer en utilisant l'hypoth\`ese de r\'ecurrence
que si
$$
\begin{array}{c}
\sigma = \sigma_0 + \mu\\
p_1 \geq p + \mu r\\
q_1 \geq q + \sgn (\mu) + (r-1)\mu,
\end{array}
$$
alors $H^{p_1,q_1}(X,\produit{\sigma}{E})=0$.

Ceci d\'ecoule du fait que le lemme \ref{E_infini_nul} reste vrai et 
que l'on a l'analogue des lemmes
\ref{ok_recurrence_q} et \ref{q_suffisant}:

\begin{lemm}
Dans l'ensemble ordonn\'e $\cal N$, on a : $(n-q_1,\sigma) < (n-q,\sigma_0)$.
\label{ordre_q}
\end{lemm}
\dem
On peut aussi supposer que $q_1=q + \sgn (\mu) + (r-1)\mu$ et 
$\sigma = \sigma_0 + \mu$. Alors si $\mu < 0$, $q_1 \geq q + r \mu$ et le lemme
\ref{ordre} s'applique. Si $\mu = 1$, ce m\^eme lemme fonctionne car $q_1 = q + r$.
Pour $\mu>1$,  en posant $x=n-(q+r)$,
ce lemme donne 
$(n-(q+r)-(\mu-1)\delta(n-(q+r)),\sigma+\mu)<(n-(q+r),\sigma + 1)$. Comme
$\delta(n-(q+r)) \leq \delta(n-q)-1$, on a
$(n-(q+r)-(\mu-1)(r-1),\sigma+\mu)<(n-(q+r),\sigma + 1)$.L'\'egalit\'e
$n-q_1= n - q - 1 - \mu (r-1) = n-(q+r)-(\mu-1)(r-1)$ donne alors
$(n-q_1,\sigma)<(n-q,\sigma_0)$, et le lemme est d\'emontr\'e.
\fin

\begin{lemm}
Si $p>P(q,\sigma_0)$ alors $p_1>P(q_1,\sigma)$.
\end{lemm}
\dem 
On peut supposer que $p_1=p+\mu r$ et $q_1=q + \mu (r-1) + \sgn (\mu)$. 
Calculons alors
$P(q,\sigma_0) - P(q + \mu (r-1) + \sgn (\mu),\sigma_0+\mu) + \mu r$. 
Cette quantit\'e vaut:
$$
\begin{array}{c}
2\mu r + 2 \mu \sigma_0 + \mu^2 + \sgn (\mu)+ (r+\sigma_0)(ae-k+2\sigma_0)\\
-[\delta(n-q-\sgn(\mu)-\mu (r-1)) + \sigma_0 + \mu][ae-k + 2\sigma_0 + 2\mu].
\end{array}
$$

\noindent
En remarquant
que $2\mu r + 2 \mu \sigma_0 =2\mu(r+\sigma_0)$ et en factorisant par
$r + \sigma_0$, on voit que cette quantit\'e \'egale
$$
\begin{array}{c}
\mu^2 + \sgn(\mu) + (ae-k+2\sigma_0+2\mu)(r+\sigma_0)\\
-(ae-k+2\sigma_0+2\mu)
[\delta(n-q-\sgn(\mu)-\mu (r-1)) + \sigma_0 + \mu],
\end{array}
$$
soit $(ae-k+2\sigma)[r-\mu-\delta(n-q-\sgn(\mu)-\mu (r-1))]+\mu^2-\sgn (\mu)$.

Comme $ae-k+\sigma \geq 0$ et que le lemme \ref{ordre_q} implique\\
$\delta(n-q-\sgn(\mu)-\mu(r-1))+\mu \leq \delta(n-q)= r$, 
ce dernier nombre est donc positif.
\fin

\para
On peut maintenant regrouper les propositions \ref{proposition_q} et 
\ref{proposition_p} sous la
forme du
\begin{theo}
Soit $E$ un fibr\'e ample de rang $e$ sur une vari\'et\'e $X$ de dimension $n$, 
$a\in \N$,
$(k_i)_{1 \leq i \leq a}$ et $(\alpha_i)_{1 \leq i \leq a}$ des entiers tels
que $\alpha_i < k_i$; soit $\sigma=\sum \alpha_i$ et $k=\sum k_i$, alors
$$ 
\begin{array}{c}
H^{p,q} (X, \otimes_i Z^{k_i-\alpha_i-1,k_i} E)=0 \\
\mbox{ si }\\
p + q > n + [\min [\delta(n-p),\delta(n-q)]+\sigma][ae-k+2\sigma] - 
\sigma(\sigma+1). 
\end{array}
$$
\label{nahm+}
\end{theo}

%******************************************************************************

%******************************************************************************

\section{R\'esultats topologiques en petit co-rang}

%******************************************************************************

\subsection{R\'esultats topologiques}

L'objet de ce paragraphe est la d\'emonstration du
\begin{theo}
Soit $E$ un fibr\'e vectoriel
de rang $e$, sur une vari\'et\'e $X$ lisse, projective et de dimension $n$.
Supposons que $E$ est
muni d'une forme quadratique \`a valeurs dans un fibr\'e en droites $L$. 
Soit $k$ un entier et
supposons que
\begin{itemize}
\item
$\dim D_k(E)=\rho:=n-t(e-k)$.
\item
$S^2E^* \otimes L$ est ample.
\item
$e-k \leq 4$.
\end{itemize}
Alors, l'application de restriction 
$H^q(X,{\cal O}_X) \rightarrow H^q(D_k(E),{\cal O}_{D_k(E)})$ est un
isomorphisme pour $0 \leq q<\rho-1$, et est injective pour $q=\rho-1$.
\label{degenerescence}
\end{theo}
\dem
Notons $D:=D_k(E)$. Tout d'abord,
il existe une r\'esolution du faisceau ${\cal O}_D$ sur ${\cal O}_X$ par des
fibr\'es vectoriels: 
le $i$-i\`eme terme $R^i$ d'une telle r\'esolution est donn\'e par
$$
R^i:=
\bigoplus_
{\begin{array}{c}
\lambda=(2l,\mu,\mu^*)\\
|\mu|+l(2l-1)=i
\end{array}}
S_{\lambda(k-1)} E \otimes L^{-l(2l+k-1)}.$$
Dans cette formule, l'expression $\lambda=(2l,\mu,\mu^*)$ signifie que la
partition $\lambda$ est de rang $2l$, que $\lambda_i=2l+\mu_i$ si 
$1 \leq i \leq 2l$ et $\lambda_{i+2l}=\mu^*_i$. Si $\lambda$ est
une partition, $\lambda(k-1)$ est obtenue en intercalant $k-1$ parts \'egales au
rang de $\lambda$ \`a $\lambda$: si par exemple $\lambda=\yng(3,2,1)$, alors
$\lambda(2)=\yng(3,2,2,2,1)$. Si $k=0$, alors $(2l,\mu,\mu^*)(-1)$ est la
partition $\nu$ telle que $\nu_i=2l+\mu_i$ pour $1 \leq i \leq 2l-1$,
$\nu_{2l}=\mu^*_1+\mu_{2l}$, et $\nu_{2l-1+i}=\mu^*_i$ pour $i \geq 2$. 
\lpara

L'existence de cette r\'esolution, bien connue des sp\'ecialistes, peut se 
justifier de la fa\c{c}on suivante: soit $Y$ l'espace total du fibr\'e
$S^2 E^* \otimes L$ et $D \subset Y$ 
le sch\'ema des formes sym\'etriques de rang au
plus $k$. Par  \cite[th\'eor\`eme 3.19, p.139]{jozefiak} et \cite{nielsen}, 
on a une r\'esolution de $\oo_D$ par des
$\oo_Y$-modules localement libres analogue \`a celle que
j'ai d\'ecrite. La section $s$ de $S^2E^* \otimes L$ induit un morphisme
$X \rightarrow Y$ et, comme l'explique Nielsen \cite{nielsen}, on peut tirer
en arri\`ere cette r\'esolution pour obtenir la r\'esolution de $\oo_{D_k(E)}$
souhait\'ee. Pour cela, il suffit en effet
d'annuler les faisceaux $\underline{Tor}^Y_i(\oo_D,\oo_X)$ pour $i>0$, ce qui
r\'esulte d'une version du corollaire 1.10 de \cite{peskine}, o\`u l'on peut
remplacer, avec leurs notations, l'hypoth\`ese 
$prof(B_{\mathfrak p}) \geq prof(A_{\varphi^{-1}({\mathfrak p})})$ par 
l'hypoth\`ese $prof(B_{\mathfrak p}) \geq dp_A(M)$. Pour annuler
$\underline{Tor}^Y_i(\oo_D,\oo_X)$, on applique
ce corollaire \`a $Spec(B)$ un ouvert 
affine de $X$ trivialisant $S^2E^* \otimes L$, $Spec(A) \subset Y$ 
l'image r\'eciproque de $Spec(B)$ par la projection
$Y \rightarrow X$ et $M$ le $A$-module d\'efinissant le faisceau $\oo_D$
sur $Spec(A)$. Pour tout id\'eal
$\mathfrak p \subset B$ correspondant \`a un point du support de
$M \otimes_A B$, on a $prof(B_{\mathfrak p}) \geq t(e-k)$, 
puisque $X$ est suppos\'ee lisse et $\dim D_k(E)=n-t(e-k)$. L'hypoth\`ese 
$prof(B_{\mathfrak p}) \geq dp_A(M)$ est donc bien v\'erifi\'ee.
\lpara

Convenons qu'une partition $\lambda$ de la forme 
$(2l,\mu,\mu^*)(k-1)$ sera dite
{\it $(k-1)$-sym\'etrique} (si $k-1= 0$, ceci signifie que $\lambda^*=\lambda)$.
Si $\lambda$ est $(k-1)$-sym\'etrique, notons $i(\lambda,k)$ l'entier tel que
$S_{\lambda} E^* \otimes L^{-l(2l+k-1)} \subset R^{i(\lambda,k)}$.

\lpara

En utilisant cette r\'esolution, il est facile de montrer que le th\'eor\`eme
\ref{degenerescence} est cons\'equence de la proposition suivante
(cet argument est par exemple expliqu\'e dans \cite{manivel}):
\begin{prop}
Soit $F$ un fibr\'e vectoriel ample de rang $e$
sur une vari\'et\'e $X$ projective lisse
de dimension $n$. Soit $k$ un entier tel que $0 < e-k \leq 4$
et $\lambda$ une partition 
$(k-1)$-sym\'etrique. Alors,
$H^{n,q}(X,S_\lambda F)=0$ si $q>t(e-k)+1-i(\lambda,k)$.
\label{petite_partition}
\end{prop}
\dem
Puisque $e-k \leq 4$, le rang $2l$ d'une partition $(k-1)$-sym\'etrique et de
longueur inf\'erieure ou \'egale \`a $e$ est n\'ecessairement inf\'erieur ou \'egal \`a 4;
on a donc $l=1$ ou $l=2$. Les partitions obtenues avec $l=2$ 
ne nous posent pas de probl\`eme, car le
th\'eor\`eme A' de \cite{manivel} donne la borne $10-i$. Il suffit donc de traiter
le cas $l=1$.

Pour simplifier les notations, j'ai suppos\'e, dans le tableau suivant, 
que de plus $k=1$, et l'ai list\'e toutes les partitions 
$\lambda$ \`a consid\'erer, et indiqu\'e, en
face, l'entier $q_0$ tel que si $q>q_0$, alors $H^{n,q}(X,S_\lambda E)=0$
pour tout fibr\'e $E$ ample
de rangs 4 et 5. La derni\`ere colonne indique le num\'ero du lemme qui montre
l'annulation du groupe de cohomologie pour la borne indiqu\'ee. L'indication
A' renvoie au th\'eor\`eme A' de
\cite{manivel}.
La colonne interm\'ediaire indique la valeur
de $q_0$ maximale pour montrer le th\'eor\`eme \ref{degenerescence}, qui vaut,
comme nous venons de l'\'etablir,
$t(e-k)+1-i$. Remarquons que lorsque $k$ prend d'autres valeurs que 1 mais
que $e-k$ reste \'egal \`a 3 ou 4, puisque
les lemmes indiqu\'es donnent une borne qui ne d\'epend que de $e-k$, on
obtient encore le r\'esultat de la proposition \ref{petite_partition}.

$$
\begin{array}{cc}
\label{tableau}
rang(E)=4& 
rang(E)=5\\
\\
\hspace{-.2cm}
\begin{tabular}{|c||c|c|c|}
\hline
\mbox{partition} & {$q_0$} & 7-i & \mbox{lemme}\\
\hline
\parbox{2.5cm}{\yng(2,2) \vspace{.5mm}} & 4 & 6 & A'\\
\hline
\parbox{2.5cm}{\yng(3,2,1) \vspace{.5mm}} & 
\parbox{.2cm}{5\\4} & 5 & \parbox{1cm}{\ref{l1}\\\ref{l1+}}\\
\hline
\parbox{2.5cm}{\yng(3,3,2) \vspace{.5mm}} & 4 & 4 & \ref{l3}\\
\hline
\parbox{2.5cm}{\yng(4,2,1,1) \vspace{.5mm}} & 2 & 4 & A'\\
\hline
\parbox{2.5cm}{\yng(4,4,2,2) \vspace{.5mm}} & 0 & 2 & A'\\
\hline
\parbox{2.5cm}{\yng(4,3,2,1) \vspace{.5mm}} & 2 & 3 & \ref{l5}\\
\hline
\end{tabular}
&
\begin{tabular}{|c||c|c|c|}
\hline
\mbox{partition} & {$q_0$} & 11-i & \mbox{lemme}\\
\hline
\parbox{2.6cm}{\yng(5,5,2,2,2) \vspace{.5mm}}& 0 & 4 & A'\\
\hline
\parbox{2.6cm}{\yng(5,2,1,1,1) \vspace{.5mm}} & 3 & 7 & A'\\
\hline
\parbox{2.6cm}{\yng(2,2) \vspace{.5mm}} & 6 & 10 & A'\\
\hline
\parbox{2.6cm}{\yng(3,2,1) \vspace{.5mm}} & 8 & 9 & \ref{l1}\\
\hline
\parbox{2.6cm}{\yng(5,4,2,2,1) \vspace{.5mm}} & 
\parbox{.2cm}{4 \\3} & 5 & 
\parbox{1cm}{\ref{l2} \\ \ref{l2+}}\\
\hline
\parbox{2.6cm}{\yng(3,3,2) \vspace{.5mm}} & 8 & 8 & \ref{l3}\\
\hline
\parbox{2.6cm}{\yng(4,2,1,1) \vspace{.5mm}} & 8 & 8 & \ref{l3}\\
\hline
\parbox{2.6cm}{\yng(4,4,2,2) \vspace{.5mm}} & 6 & 6 & \ref{l4}\\
\hline
\parbox{2.6cm}{\yng(5,3,2,1,1) \vspace{.5mm}} & 4 & 6 & \ref{l5}\\
\hline
\parbox{2.6cm}{\yng(4,3,2,1) \vspace{.5mm}} &
\parbox{.2cm}{7 \\6} & 7 &
\parbox{1.5cm}{\ref{l6} \\ \ref{l7}}\\
\hline
\end{tabular}
\end{array}
$$

%******************************************************************************

\newpage
\subsection{Th\'eor\`emes d'annulation}

Dans ce paragraphe, je montre les lemmes annonc\'es dans le pr\'ec\'edent. Soit $X$
une vari\'et\'e projective lisse de dimension $n$,
$E$ un fibr\'e ample de rang $e$ sur $X$, et $L$ un fibr\'e en droites
nef.

Pour des raisons qui sont 
expliqu\'ees apr\`es
le lemme \ref{l2}, je note $\lambda=[a,b,c,d]$ la partition de rang
2 telle que $\lambda_1=a+2,\lambda_2=b+2,\lambda^*_1=e-c,\lambda_2^*=e-d$.
Par exemple, $S_{[0,0,0,0]} E=(\det E)^2$. De même, je note $\lambda=[a,b]$
la partition de rang 1 telle que $\lambda_1=a+1$ et $\lambda_1^*=e-b$. Par
la suite, $a\leq b,c\leq d,t,q$ sont des entiers positifs ou nuls.

\begin{lemm}
\ \vspace{-.56cm}
$$
\begin{array}{rcl}
\hspace{2cm}H^{n,q}(X,S_{[1,c+1]}E \otimes S_{[0,d]}E \otimes L^t) & =0 & 
\mbox{ si } q>2c+d+2, \mbox{ et }\\
H^{n-1,q}(X,S_{[0,c]}E \otimes S_{[0,d]}E \otimes L^t) &=0 & 
\mbox{ si } q>2c+d+1.
\end{array}
$$
\label{l1}
\end{lemm}
\ \vspace{-1cm}\\
\rem
Le th\'eor\`eme \ref{nahm+} donne la borne $2c+2d+2$ au lieu de $2c+d+2$
et les autres r\'esultats que je
connais ont une borne qui d\'epend de $k$ (par exemple, 
\cite[th\'eor\`eme A]{manivel} donne la borne $k+c+d$); il est clair que pour 
d\'emontrer le th\'eor\`eme \ref{degenerescence} par cette m\'ethode, il faut une borne
ind\'ependante de $k$.\\
\dem
Pour montrer ce lemme, je vais appliquer la technique, 
nouvelle \`a ma connaissance,
de ``comparaison de suites spectrales''.

Tout d'abord, on peut supposer que $e-c$ est pair. En effet, supposons
que ce lemme est
vrai pour $e-c$ pair. Si $e-c$ est impair
et $M$ est un fibr\'e en droites ample sur $X$, alors on peut consid\'erer
le fibr\'e vectoriel $E \oplus M$ de rang $e+1$; comme $e+1-c$ est alors pair,
on peut appliquer le lemme, et comme pour le lemme \ref{astuce_nahm}, on en
d\'eduit que $H^{n,q}(X,S_{[1,c+1]}E \otimes S_{[0,d]}E \otimes L^t)=0$ 
si $q>2c+d+2$ et $H^{n-1,q}(X,S_{[0,c]}E \otimes S_{[0,d]}E \otimes L^t)=0$
si $q>2c+d+1$.

Soit alors un entier $l$ tel que $e-c=2l$. Consid\'erons
d'abord la vari\'et\'e $Y=\p E^* \times_X \p E^*$, notons $\pi:Y \rightarrow X$ la
projection naturelle, et consid\'erons le fibr\'e en droites 
${\cal L}={\cal O}(2l,e-d) \otimes \pi^* L^t$ sur $Y$. Comme
${\cal O}(2l,e-d)$ est ample (en effet, $E$ est ample) 
et $\pi^*L$ nef, $\cal L$
est ample.
Soit aussi $P=n+(2l-1)+(e-d-1)-1$. 

D\'eterminons les groupes $\Ed{P}{t}{i}{j}{2}$ de la suite spectrale de Leray
introduite au paragraphe \ref{para_leray}. Ceux-ci valent, par la formule
(\ref{leray}), $H^{t,i}[X,R^{P-t,j}\pi_*\ooo(2l,e-d) \otimes L^t]$.
Ces groupes sont donc nuls si $t>n$ ou $P-t>(2l-1)+(e-d-1)$. En effet, si $V$
est un espace vectoriel fix\'e et $l$ et $m$ sont des entiers positifs,
alors $H^{i,j}[\p V,\ooo(l)]=0$ si
$i \geq l$ (cons\'equence directe de la proposition 
\ref{partition_admissible}), et donc 
$H^{i,j}[\p V \times \p V,\ooo(l,m)]=0$ si $i \geq l+m$ 
(formule de K\" unneth). Les seules valeurs possibles sont donc
$t=n$ ou $t=n-1$, et on a
$$
\begin{array}{rcl}
\Ed{P}{n-1}{i}{0}{2} & = & 
H^{n-1,i}(X,S_{[0,c]}E \otimes S_{[0,d]}E \otimes L^t),
\mbox{\ \ \  et }\\
\Ed{P}{n}{i}{0}{2} & = & 
H^{n,i}(X,L^t \otimes R^{(2l-1)+(e-d-1)-1,0}\pi_* \ooo(2l,e-d))\\
& = & H^{n,i}(X,S_{[1,c+1]}E \otimes S_{[0,d]}E \otimes L^t)\\
& \oplus & H^{n,i}(X,S_{[0,c]}E \otimes S_{[1,d+1]}E \otimes L^t).
\end{array}
$$

\noindent
Ces suites spectrales sont donc d\'eg\'en\'er\'ees, et on en d\'eduit que
$$
\begin{array}{rcl}
\Es{P}{n-1}{q-n+1}{1} & = &H^{n-1,q}(X,S_{[0,c]}E 
\otimes S_{[0,d]}E \otimes L^t),\\
\Es{P}{n}{q-n}{1} & = & H^{n,q}
(X,S_{[1,c+1]}E \otimes S_{[0,d]}E \otimes L^t)\\
& \oplus & H^{n,q}(X,S_{[0,c]}E \otimes S_{[1,d+1]}E \otimes L^t),
\end{array}
$$
et les autres
termes sont nuls. Par ailleurs,
$\Es{P}{p}{q-p}{\infty}$ est une composante de
$H^{P,q}(X,{\cal L})$; il est donc nul par le th\'eor\`eme de Kodaira
si 
$$P+q>n+2(e-1), \mbox{ soit } q>n+2(e-1)-P=c+d+1.$$
On en d\'eduit donc que dans ce
cas, la flêche connectant $\Es{P}{n-1}{q-n+1}{1}$ et
$\Es{P}{n}{q-n+1}{1}$ est un isomorphisme, ce qui conduit \`a:

$$
\begin{array}{rl}
& H^{n-1,q}(X,S_{[0,c]}E \otimes S_{[0,d]} E \otimes L^t)\\ 
= & H^{n,q+1}(X,S_{[1,c+1]}E \otimes S_{[0,d]}E \otimes L^t)\\
\oplus & H^{n,q+1}(X,S_{[0,c]}E \otimes S_{[1,d+1]}E \otimes L^t)
\mbox{ si } q>c+d+1.
\end{array}
$$

\lpara

En consid\'erant une autre suite spectrale, on va obtenir une autre \'egalit\'e
faisant intervenir ces termes; en comparant ces \'egalit\'es, on prouvera qu'ils
sont nuls. Soit donc maintenant $Y=G_2(E^*) \times_X \p E^*$ et
${\cal L}={\cal O}(l,e-d)$ sur $Y$, et posons $P=n+3(l-1)+(e-d-1)-1$. 
On a alors
de fa\c{c}on similaire
$\Es{P}{n-1}{q-n+l}{1}=H^{n-1,q}(X,S_{[0,c]}E \otimes S_{[0,d]}E \otimes L^t)$,
mais par contre,
$\Es{P}{n}{q-n+l-1}{1}=H^{n,q}(X,S_{[0,c]}E \otimes S_{[1,d+1]}E \otimes L^t)$.
Ceci conduit \`a l'\'egalit\'e
$$
H^{n-1,q}(X,S_{[0,c]}E \otimes S_{[0,d]}E \otimes L^t)=
H^{n,q+1}(X,S_{[0,c]}E \otimes S_{[1,d+1]} E \otimes L^t)
\mbox{ si } q>2c+d+1.
$$

\lpara

En comparant ces \'egalit\'es, on obtient
$$H^{n,q}(X,S_{[1,c+1]}E \otimes S_{[0,d]} E \otimes L^t)=0,$$
soit la premi\`ere affirmation du lemme.
En utilisant la r\`egle de \litt, on va montrer l'autre r\'esultat du lemme.
En effet, rappelons que $c \leq d$; si $c=d$, on peut conclure; 
supposons donc $c<d$. Cette r\`egle implique que
$$
\begin{array}{rcl}
S_{[1,c+1]}E \otimes S_{[0,d]}E & = & 
\bigoplus_{0 \leq x \leq c+1} S_{[1,0,x,d+c+1-x]}E \bigoplus R_1\\
S_{[0,c]}E \otimes S_{[1,d+1]}E & = & 
\bigoplus_{0 \leq x \leq c} S_{[1,0,x,c+d+1-x]}E \bigoplus R_2,\\
\end{array}
$$
o\`u $R_1$ et $R_2$ sont des sommes de composantes de type
$S_{[0,0,x,y]}$, que l'on sait d\'ej\`a annuler. L'annulation de la cohomologie
de $S_{[1,c+1]}E \otimes S_{[0,d]}E \otimes L^t$ entraine donc celle de
$S_{[0,c]}E \otimes S_{[1,d+1]}E \otimes L^t$, 
et le lemme est alors cons\'equence d'une
des \'egalit\'es d\'emontr\'ees.
\fin

\para

Notons que ce lemme montre par exemple
$H^{n,q}(X,S_{[1,0,c,c+1]}E)=0$ si $q>3c+2$, car
$S_{[1,0,c,c+1]}E \subset S_{[1,c+1]}E \otimes S_{[0,c]}E$. Ce n'est pas la
borne indiqu\'ee dans le tableau p. \pageref{tableau}.
En effet, pour cette partition particuli\`ere, 
on peut un peu raffiner le raisonnement pr\'ec\'edent. Auparavant, je souhaite
faire une remarque d'ordre g\'en\'eral qui all\`egera les calculs.

En g\'en\'eralisant le raisonnement utilis\'e pour le lemme pr\'ec\'edent,
on voit que l'on obtient simultan\'ement l'annulation de
$H^{p_1,q}(X,S_{\lambda_1} E)$ et celle de
$H^{p_2,q}(X,S_{\lambda_2} E)$; supposons que $p_1<p_2$: la borne $q_1$
obtenue pour le groupe $H^{p_1,q}(X,S_{\lambda_1} E)$ est la
diff\'erence entre la dimension de $Y$ et la somme $i+j+p_1$, o\`u
$i$ et $j$ sont les entiers tels que
si $V$ est un espace vectoriel fix\'e, $Y_V$ la fibre
de la projection $Y \rightarrow X$, et ${\cal L}_V$ la restriction du fibr\'e en
droites $\cal L$ sur une telle fibre, alors 
$H^{i,j}(Y_V,{\cal L}_V)=S_\lambda V$. La borne $q_2$ pour
$H^{p_2,q}(X,S_{\lambda_2} E)$ vaut $q_1+t$, o\`u $t$ est l'entier tel
que les flêches de la suite spectrale connectent
$H^{p_1,q}(X,S_{\lambda_1} E)$ et
$H^{p_2,q+t}(X,S_{\lambda_2} E)$. Attention,
cette recette n'est bien sûr valable
que si aucun autre terme de la suite spectrale ne vient compenser les termes
que l'on veut annuler. Par ailleurs, on peut facilement exprimer cette
diff\'erence:

\begin{rema}
Une suite spectrale donn\'ee par un fibr\'e en droites sur
$G_r(E^*) \times_X G_s(E^*)$ peut montrer l'annulation de
$H^{p,q}(X,S_{[a,b]}E \otimes S_{[c,d]}E \otimes L^t)$ pour
$q>(r-1)a+rb+(s-1)c+sd$.
\label{rema_borne}
\end{rema}
\dem
En effet, comme je l'ai expliqu\'e, il s'agit d'\'evaluer la diff\'erence entre
$\dim G_r(V)$ et $i+j$, o\`u $i$ et $j$ sont les entiers tels que
$H^{i,j}(G_r(V),{\cal O}(l))=S_{[a,b]}V$. Cette diff\'erence vaut
$ra+(r-1)b$, car $i=(l-1)t(r)-ar,j=(l-1)t(r-1)-a(r-1)$, et
$\dim V=rl-a+b$.
\fin

Dans le lemme qui suit, je traite des partitions qui sont un cas particulier de
celles trait\'ees par le lemme \ref{l1}; pour ces partitions particuli\`eres, je
peux r\'eduire de 1 la borne donn\'ee par le lemme \ref{l1}.

\begin{lemm}
Supposons $c>0$. \vspace{-.1cm}
$$
\begin{array}{rcl}
H^{n,q}(X,S_{[1,0,c+1,c]}E \otimes L^t) & =0 & 
\mbox{ si } q>3c+1, \mbox{ et }\\
H^{n-1,q}(X,S_{[0,0,c,c]}E \otimes L^t) &=0 & 
\mbox{ si } q>3c.
\end{array}
$$
\label{l1+}
\end{lemm}
\ \vspace{-.8cm}\\
\dem
Tout d'abord, on peut supposer que $e-c$ est impair et
$e-c>>0$. Dans la d\'emonstration du lemme \ref{l1}, 
j'ai montr\'e comment tenir compte
de $L^t$; pour all\'eger les notations, je suppose dor\'enavant que $t=0$.

Soit alors $l$ un entier tel que $e-c=2l+1$. Consid\'erons
d'abord la vari\'et\'e $Y=\p E^* \times_X \p E^*$, et le fibr\'e
${\cal O}(2l+1,2l+1)$ sur $Y$.
Comme dans la d\'emonstration du lemme pr\'ec\'edent, on montre que
$$
2H^{n,q+1}(X,S_{[1,c+1]}E \otimes S_{[0,c]})=
H^{n-1,q}(X,S_{[0,c]}E \otimes S_{[0,c]}E) \mbox{ si } q>2c+1.$$
Soit maintenant 
${\cal O}(2l,l+1) \rightarrow \p E^* \times_X G_2(E^*)$. On a alors de fa\c{c}on
similaire
$$H^{n,q+1}(X,S_{[1,c+2]}E \otimes S_{[0,c-1]} E)=
H^{n-1,q}(X,S_{[0,c+1]}E \otimes S_{[0,c-1]} E) 
\mbox{ si } q>3c.$$
Par la r\`egle de \litt,
$S_{[1,c+1]} E \otimes S_{[0,c]}E=S_{[1,c+2]}E \otimes S_{[0,c-1]} E
\oplus S_{[1,0,c,c+1]}E$
(\`a des termes de type $S_{[0,0,x,y]}E$ pr\`es), et
$S_{[0,c]}E \otimes S_{[0,c]}E=S_{[0,c-1]}E \otimes S_{[0,c+1]}E
\oplus S_{[0,0,c,c]}E$. 

Je donne maintenant une derni\`ere \'egalit\'e faisant
intervenir $H^{n,q+1}(X,S_{[1,0,c,c+1]}E)$
et $H^{n-1,q}(S_{[0,0,c,c]}E)$, qui permettra de conclure.
Consid\'erons ${\cal O}(e-c)\rightarrow G_2(E^*)$. Par la proposition
\ref{partition_admissible}, une partition $\lambda$, de taille au plus
$(e-2) \times 2$, $(e-c)$-admissible, v\'erifie soit 
$\lambda_1,\lambda_2 \leq e-c-2$, soit $\lambda_1=\lambda_2+e-c-1$. Dans le
deuxi\`eme cas, $|\lambda|=\lambda_1+\lambda_2\leq [e-2-(e-c-1)]+[e-2]=e+c-3$.
Par contre, la partition $(e-c-2,e-c-2)$ est de poids $2(e-c-2)$. Il est clair
qu'\`a $c$ fix\'e, si $e$ est grand, le poids de cette partition est strictement
sup\'erieur \`a toute constante plus le 
poids d'une partition du deuxi\`eme type. Quitte \`a
augmenter $e$, on peut donc faire comme si seules les partitions du premier
type intervenaient. En posant $P=n+2(e-c-2)-1$, on montre alors que
$$
H^{n,q+1}(X,S_{[1,0,c,c+1]}E)=H^{n-1,q}(X,S_{[0,0,c,c]}E) \mbox{ si }
q>2c+1.
$$

Les trois \'egalit\'es montrent que si $q>3c$, alors 
$H^{n,q+1}(X,S_{[1,c+1]}E \otimes
S_{[0,c]}E)=0$. La premi\`ere implique alors le reste du lemme.
\fin

\begin{lemm}
\ \vspace{-.52cm}
$$
\begin{array}{rcccl}
\parbox{1.5cm}{\ }&H^{n,q}(X,S_{[a,0]}E \otimes S_{[b+1,1]}E \otimes L^t) & 
= & 0 & \mbox{ si } q>a+2, \mbox{ et }\\
&H^{n-1,q}(X,S_{[a,0]} \otimes S_{[b,0]}E \otimes L^t)& = & 0 &
\mbox{ si } q>a+1.
\label{l2}
\end{array}
$$
\end{lemm}
\noindent
\rem
Ce lemme me donne l'occasion de justifier
mes notations. Les partitions les plus
faciles \`a traiter sont celles qui correspondent \`a une puissance du d\'eterminant:
en effet, puisque $(\det E)^l$ est ample pour tout $l$, le th\'eor\`eme de
Kodaira donne directement que
$H^{p,q}(X,(\det E)^l)=0$ si $p+q>n$. Lorsqu'on s'\'ecarte de cette partition,
je vais expliquer qu'il existe
une sym\'etrie entre le fait d'allonger les premi\`eres parts et le fait de
creuser le bas de la partition. Autrement dit, 
dans la notation $\lambda=[a,b,c,d]$, $a$ et $c$, et $b$ et $d$, jouent des
rôles sym\'etriques. On peut d'ailleurs remarquer que le th\'eor\`eme 
\ref{nahm+} n'\'echappe pas \`a ce principe: si l'on note 
$r=\min \{\delta(n-p),\delta(n-q)\}$, alors la borne vaut
$q_0=n-p+(r+\sigma)(ae-k+2\sigma)-\sigma(\sigma+1)$. 
Or, avec mes notations, on a
$Z^{k_i-\alpha_i-1,k_i}=[\alpha_i,e-k_i+\alpha_i]$. Notons 
$\beta_i=e-k_i+\alpha_i$; on a $ae-k+\sigma=\sum_i (e-k_i+\alpha_i)=
\sum_i \beta_i=:\tau$. Ainsi, on a donc
\begin{equation}
q_0=n-p+(r+\sigma)(\tau+\sigma)-\sigma(\sigma+1)
=n-p+r(\sigma+\tau)+\sigma(\tau-1).
\label{nahm_symetrique}
\end{equation}
On voit bien que cette formule est sym\'etrique en $\sigma$ et $\tau\ (-1)$. Le
lecteur intrigu\'e pourra s'amuser de constater que la d\'emonstration
du th\'eor\`eme \ref{nahm+} fonctionne de fa\c{c}on tout \`a fait similaire si l'on
fait jouer par $\tau$ le rôle jou\'e par $\sigma$, 
et donne le même r\'esultat. Il est vraisemblable que
le d\'emonstration de th\'eor\`emes d'annulation efficaces pour des partitions de
rang strictement sup\'erieur \`a 1 utilise une r\'ecurrence qui fasse intervenir de
fa\c{c}on combin\'ee $\sigma$ et $\tau$. 

Le lemme \ref{l2}, pour cette sym\'etrie, est le
sym\'etrique du lemme \ref{l1}; il se d\'emontre de fa\c{c}on tout \`a fait similaire.\\
\dem
On se ram\`ene comme pr\'ec\'edemment au cas $e+a$ pair, puis
on consid\`ere
le fibr\'e en droites ${\cal O}(e+a,e+b)$ sur $\p E^* \times_X \p E^*$,
et le fibr\'e en droites ${\cal O}(e+b,\frac{e+a}{2})$ sur
$\p E^* \times_X G_2(E^*)$. On utilise aussi, comme pour d\'emontrer
le lemme \ref{l1}, la r\`egle de \litt. En tenant compte de la remarque
\ref{rema_borne}, on trouve ainsi
$H^{n-1,q}(X,S_{[a,0]} \otimes S_{[a+1,1]}E)=0$ si 
$q>a+1$.
Dans la
suite spectrale, il y a une fl\`eche entre
$H^{n,q+1}(X,S_{[a,0]}E \otimes S_{[a+1,1]}E)$ et
$H^{n-1,q}(X,S_{[a,0]} \otimes S_{[a+1,1]}E)$; 
$H^{n,q}(X,S_{[a,0]}E \otimes S_{[a+1,1]}E)$ s'annule donc si
$q>a+2$.
\fin

Par ailleurs, on peut, comme le fait le lemme \ref{l1+}, raffiner un peu le
lemme pr\'ec\'edent pour traiter certains cas particuliers:
\begin{lemm}
Supposons $a>0$. \vspace{-.1cm}
$$
\begin{array}{rcl}
H^{n,q}(X,S_{[a+1,a,0,1]}E \otimes L^t) & =0 & 
\mbox{ si } q>a+1, \mbox{ et }\\
H^{n-1,q}(X,S_{[a,a,0,0]}E \otimes L^t) & =0 & 
\mbox{ si } q>a.
\end{array}
$$
\label{l2+}
\end{lemm}
\dem
Similaire \`a celle du lemme \ref{l1+}.
\fin

\lpara
\begin{lemm}\hspace{-.2cm}
$H^{n,q}(X,S_{[1,c+1]}E \otimes S_{[1,c+1]}E \otimes L^t)=
H^{n,q}(X,S_{[0,c+2]}E \otimes S_{[0,c]}E \otimes L^t)=0$
si $q>4c+4$, et
$H^{n-1,q}(X,S_{[1,c+1]} \otimes S_{[0,c]}E \otimes L^t) = 0$
si $q>4c+3$.
\label{l3}
\end{lemm}
\dem
Comme pr\'ec\'edemment, on peut supposer que $e-c$ est pair et $t=0$; soit donc
$l$ tel que $e-c=2l$. Consid\'erons le fibr\'e en droites ${\cal O}(2l,2l)$ sur
$\p E^* \times_X \p E^*$, et la suite spectrale 
correspondante pour $P=n+2(2l-1)-2$.
Les termes $\Es{P}{i}{j}{1}$ ne peuvent 
être non nuls que si $i$ vaut $n,n-1$, ou
$n-2$. Or $\Es{P}{n-2}{q-n+2}{1}=H^{n-2,q}(X,S_{[0,c]} E 
\otimes S_{[0,c]} E)$, qui est nul si $q>4c+2$ par le th\'eor\`eme
\ref{nahm+}. On en d\'eduit donc que si $q>4c+2$, alors
$$
\begin{array}{ll}
 & H^{n,q+2}(X,S_{[1,c+1]}E \otimes S_{[1,c+1]}E) \oplus
2H^{n,q+2}(X,S_{[2,c+2]}E \otimes S_{[0,c]}E)\\
= & 2H^{n-1,q+1}(X,S_{[1,c+1]} \otimes S_{[0,c]}E).
\end{array}
$$
De même, en consid\'erant 
${\cal O}(2l,l)\rightarrow \p E^* \times_X G_2(E^*)$, on
montre que 
\cent{$H^{n,q+2}(X,S_{[2,c+2]}E \otimes S_{[0,c]}E)=
H^{n-1,q+1}(X,S_{[1,c+1]} \otimes S_{[0,c]}E)$ si $q>4c+2$.}
Ce syst\`eme donne
$H^{n,q+2}(X,S_{[1,c+1]}E \otimes S_{[1,c+1]}E)=0$. Mais comme
$$S_{[1,c+1]}E \otimes S_{[1,c+1]}E 
\supset S_{[2,c+2]}E \otimes S_{[0,c]}E,$$
on en d\'eduit qu'aussi
$H^{n,q+2}(X,S_{[2,c+2]}E \otimes S_{[0,c]}E)=0$, et donc\\
$H^{n-1,q+1}(X,S_{[1,c+1]} \otimes S_{[0,c]}E)=0$.
\fin

\lpara

Le sym\'etrique du lemme pr\'ec\'edent est:
\begin{lemm}
$H^{n,q}(X,S_{[a+1,1]}E \otimes S_{[a+1,1]}E \otimes L^t)=
H^{n,q}(X,S_{[a+2,2]}E \otimes S_{[a,0]}E \otimes L^t)=0$ si $q>2a+4$ et
$H^{n-1,q}(X,S_{[a+1,1]} \otimes S_{[a,0]}E \otimes L^t)=0$ si $q>2a+3$.
\label{l4}
\end{lemm}
\dem
Comme pour le lemme pr\'ec\'edent, on pose $2l=a+e$ et $t=0$,
et on compare les suites spectrales correspondant
\`a ${\cal O}(2l,2l) \rightarrow \p E^* \times_X \p E^*$ et
${\cal O}(2l,l) \rightarrow \p E^* \times_X G_2(E^*)$.
\fin

Jusqu'\`a maintenant, nous avons consid\'er\'e des partitions de type
$[a,\epsilon] \otimes [b,\epsilon']$ ($\epsilon,\epsilon'=0$ ou 1), ou de type
$[\epsilon,a] \otimes [\epsilon',b]$.
Etudions finalement des partitions ``panach\'ees'', c'est-\`a-dire
de type $[a,\epsilon,b,\epsilon']$.
\begin{lemm}
Si $a$ et $c$ ont même parit\'e, alors
$H^{n,q}(X,S_{[a,0]}E \otimes S_{[1,c+1]}E \otimes L^t)=0$ si $q>2c+2$, et
$H^{n,q}(X,S_{[a+1,1]}E \otimes S_{[0,c]}E \otimes L^t)=0$ si $q>a+2$.
Si $q>\max(a,2c)+1$, alors 
$H^{n-1,q}(X,S_{[a,0]}E \otimes S_{[0,c]}E \otimes L^t)=0$.
\label{l5}
\end{lemm}
\rem
Ce lemme illustre la difficult\'e du probl\`eme qui consiste \`a obtenir des
th\'eor\`emes d'annulation optimaux: le lemme \ref{l4} avec $a=1$
donne la borne 6 pour
la partition $\yng(5,3,2,1,1)$ (si $e=5$); 
cette borne convient pour d\'emontrer le
th\'eor\`eme \ref{degenerescence} mais n'est pas optimale, puisque le lemme
\ref{l5}
donne la borne 4 ($a=3,c=1$).
Par contre, j'aurai besoin de cette bonne borne pour
montrer le lemme \ref{l6}.
Pour d\'emontrer un analogue du th\'eor\`eme \ref{degenerescence}
en tous rangs, il faut probablement prendre en compte ce genre de
subtilit\'es.\\
\dem
On suppose que $a+e=2l,e-c=2m$ et $t=0$. Si l'on consid\`ere
${\cal O}(2l,2m) \rightarrow \p E^* \times_X \p E^*$, alors on obtient
\cent{$H^{n-1,q}(X,S_{[a,0]}E \otimes S_{[0,c]}E)
=H^{n,q+1}(X,S_{[a,0]}E \otimes S_{[1,c+1]}E) \oplus
H^{n,q+1}(X,S_{[a+1,1]}E \otimes S_{[0,c]}E)$ si $q>c+1$. \vspace{.2cm}}
Avec ${\cal O}(2l,m) \rightarrow \p E^* \times_X G_2(E^*)$, on a
\cent{$H^{n-1,q}(X,S_{[a,0]}E \otimes S_{[0,c]}E)=
H^{n,q+1}(X,S_{[a+1,1]}E \otimes S_{[0,c]}E)$ si $q>2c+2$.\vspace{.2cm}}
Enfin, ${\cal O}(l,2m) \rightarrow G_2(E^*) \times_X \p E^*$, donne
\cent{$H^{n-1,q}(X,S_{[a,0]}E \otimes S_{[0,c]}E)=
H^{n,q+1}(X,S_{[a,0]}E \otimes S_{[1,c+1]}E)$ si $q>a+2$.\vspace{.2cm}}

Le lemme d\'ecoule de ces trois remarques.
\fin

\begin{lemm}
Si $a$ et $c$ ont même parit\'e, alors
\cent{$H^{n-2,q}(X,S_{[a,0]}E \otimes S_{[0,c]}E \otimes L^t)=0$ si
$q>a+2c+2$.}
\label{n-2}
\end{lemm}
\dem
Si $t=0,2l=e+a$ et $2m=e-c$, on regarde la suite spectrale associ\'ee \`a
${\cal O}(l,m) \rightarrow G_2(E^*) \times_X G_2(E^*)$, pour 
$P=n+(2e-l-3)+(3(m-1))-2$. Par la remarque \ref{rema_borne},
pour $q>a+2c+2$, celle-ci compare
$H^{n-2,q}(X,S_{[a,0]}E \otimes S_{[0,c]}E)$ et
$H^{n,q+2}(X,S_{[a+1,1]}E \otimes S_{[0,c]}E) \oplus
H^{n,q+2}(X,S_{[a,0]}E \otimes S_{[1,c+1]}E)$. Or, ces deux
derniers groupes sont
nuls, grâce au lemme \ref{l5}, si $q>\max\{a,2c\}$.
\fin

\begin{lemm}
Si $a$ et $c$ ont même parit\'e, alors
\cent{$H^{n,q}(X,S_{[a+1,1]}E \otimes S_{[1,c+1]}E \otimes L^t)=0$ si
$q>a+2c+4$.}
\label{l6}
\end{lemm}
\dem
Si $t=0$, $2l=e+a$ et $2m=e-c$, alors 
${\cal O}(2l,2m) \rightarrow \p E^* \times_X \p E^*$ donne, pour $q>a+2c+2$,
$$
\begin{array}{rl}
& H^{n,q+2}(X,S_{[a+1,1]}E \otimes S_{[1,c+1]}E) \oplus
H^{n,q+2}(X,S_{[a+2,2]}E \otimes S_{[0,c]}E)\\
\oplus & H^{n,q+2}(X,S_{[0,a]}E \otimes S_{[2,c+2]}E)\\
= & H^{n-1,q+1}(X,S_{[1,a+1]}E \otimes S_{[0,c]}E) \oplus
H^{n-1,q}(X,S_{[0,a]}E \otimes S_{[1,c+1]}E).
\end{array}
$$
En effet, $H^{n-2,q}(X,S_{[a,0]} \otimes S_{[0,c]}E)=0$ si $q>a+2c+2$ par le
lemme \ref{n-2}. De même, ${\cal O}(l,2m) \rightarrow G_2(E^*) \times_X \p E^*$
et ${\cal O}(2l,m) \rightarrow \p E^* \times_X G_2(E^*)$
donnent respectivement, pour
$q>2a+c+2$,
$$
\begin{array}{rcl}
H^{n,q+2}(X,S_{[a+2,2]}E \otimes S_{[0,c]}E) & = &
H^{n-1,q+1}(X,S_{[1,a+1]}E \otimes S_{[0,c]}E) \mbox{ et }\\
H^{n,q+2}(X,S_{[0,a]}E \otimes S_{[2,c+2]}E) & = &
H^{n-1,q}(X,S_{[0,a]}E \otimes S_{[1,c+1]}E).
\end{array}
$$
\fin
\lpara

Remarquons enfin que l'on peut am\'eliorer la borne pour
$S_{[2,1]}E \otimes S_{[1,2]}E$. En effet, en g\'en\'eralisant la d\'emonstration 
du lemme pr\'ec\'edent, on voit que l'annulation\\
$H^{n-2,q}(X,S_{[1,0]}E \otimes S_{[0,1]}E)=0$ si $q>q_0$ implique
$H^{n,q}(X,S_{[2,1]}E \otimes S_{[1,2]}E)=0$ si $q>q_0+2$. Le lemme pr\'ec\'edent
utilisait la borne $q_0=5$ montr\'ee dans le lemme \ref{n-2}; je vais montrer
que $H^{n-2,q}(X,S_{[1,0]}E \otimes S_{[0,1]}E)=0$ si $q>4$ par un argument qui
ne fonctionne que pour $S_{[1,0]}E \otimes S_{[0,1]}E$. 
La remarque sp\'ecifique est que, par
la r\`egle de \litt,
$$S_{[1,0]}E \otimes S_{[0,1]}E = 
(\det E)^2 \oplus \det E \otimes S_{[1,1]}E.$$
Puisque $H^{n-2,q}(X,(\det E)^2)=0$ si $q>2$, il suffit de montrer que
$H^{n-2,q}(X,\det E \otimes S_{[1,1]}E)=0$ si $q>4$. Pour cela, on peut
supposer que le rang de $E$ est multiple de 4, $e=4f$, et on regarde
la suite spectrale correspondant \`a
${\cal O}(f,4f) \rightarrow G_4(E) \times_X \p E$. Dans cette suite spectrale,
notre groupe, $H^{n-2,q}(X,\det E \otimes S_{[1,1]}E)=0$, est connect\'e \`a
$H^{n-3,q-1}(X,(\det E)^2)$ (qui s'annule si $q-1>3$),
$H^{n-1,q+1}(X,\det E \otimes S_{[2,2]}E)$ (nul pour $q+1>5$: lemme \ref{l5}),
et enfin $H^{n,q+1}(X,\det E \otimes S_{[3,3]}E)=0$ (nul pour $q+1>3$ par le
th\'eor\`eme A' de \cite{manivel}). On a donc d\'emontr\'e

\begin{lemm}
$H^{n,q}(X,S_{[2,1]}E \otimes S_{[1,2]}E)=0$ si $q>6$.
\label{l7}
\end{lemm}
Ceci illustre \`a nouveau la complexit\'e du sujet: on obtient de meilleures bornes
pour des partitions sp\'ecifiques que celles obtenues par la m\'ethode g\'en\'erale.
Par ailleurs, comme nous utilisons une r\'ecurrence, le fait de ne pas obtenir
le th\'eor\`eme d'annulation optimal pour une partition pr\'ecise se r\'epercute sur de
nombreuses autres partitions, entrainant progressivement la ``catastrophe''.

\end{document}